\documentclass[a4paper, fleqn]{amsart}
\usepackage{amsmath,amsfonts,amssymb,amsthm}
\usepackage{mathtools}
\usepackage{mathrsfs}
\usepackage{pdfpages}
\usepackage{url}
\usepackage{color,xcolor}
\usepackage{caption}
\usepackage{enumerate}
\usepackage{afterpage}
\usepackage{xargs,twoopt}
\usepackage{hyperref}
\usepackage[toc,page]{appendix}
\usepackage{spverbatim}
\hypersetup{colorlinks = true, linkbordercolor = {white}}
\usepackage{aliascnt}
\usepackage{bbold}
\usepackage{dsfont}
\usepackage[cal=boondoxo]{mathalfa}
\usepackage{bussproofs}
\usepackage{graphicx}
\usepackage[retainorgcmds]{IEEEtrantools}
\usepackage{stmaryrd}
\usepackage{tensor}
\usepackage{textcomp}
\usepackage{tikz}
\usepackage{units}
\usepackage[all]{xy}
\DeclareSymbolFont{extraitalic}      {U}{zavm}{m}{it}
\DeclareMathSymbol{\Qoppa}{\mathord}{extraitalic}{161}
\DeclareMathSymbol{\qoppa}{\mathord}{extraitalic}{162}
\DeclareMathSymbol{\Stigma}{\mathord}{extraitalic}{167}
\DeclareMathSymbol{\Sampi}{\mathord}{extraitalic}{165}
\DeclareMathSymbol{\sampi}{\mathord}{extraitalic}{166}
\DeclareMathSymbol{\stigma}{\mathord}{extraitalic}{168}

\newcommand{\id}{\text{id}}
\newcommand{\ceA}{\mathcal{A}}
\newcommand{\ceB}{\mathcal{B}}
\newcommand{\ceC}{\mathcal{C}}
\newcommand{\ceD}{\mathcal{D}}

\newcommand{\ceI}{\mathcal{I}}
\newcommand{\ceJ}{\mathcal{J}}

\newcommand{\ceL}{\mathcal{L}}
\newcommand{\ceM}{\mathcal{M}}

\newcommand{\ceR}{\mathcal{R}}

\newcommand{\cB}{\mathscr{B}}
\newcommand{\cC}{\mathscr{C}}
\newcommand{\cD}{\mathscr{D}}

\newcommand{\cH}{\mathscr{H}}
\newcommand{\cI}{\mathscr{I}}
\newcommand{\cJ}{\mathscr{J}}
\newcommand{\cK}{\mathscr{K}}
\newcommand{\cL}{\mathscr{L}}

\newcommand{\cR}{\mathscr{R}}

\newcommand{\LA}{\mathfrak{A}}

\newcommand{\cCat}{\mathbf{Cat}}

\newcommand{\mbC}{\mathbf{C}}

\newcommand{\mbG}{\mathbf{G}}

\newcommand{\mbK}{\mathbf{K}}

\newcommand{\mbM}{\mathbf{M}}
\newcommand{\mbN}{\mathbf{N}}

\newcommand{\bbon}{\mathbb{1}}
\newcommand{\bbk}{\mathbb{k}}

\newcommand{\bbZ}{\mathbb{Z}}

\newcommand{\tti}{\mathtt{i}}
\newcommand{\ttj}{\mathtt{j}}
\newcommand{\ttk}{\mathtt{k}}

\newcommand{\ul}[1]{\underline{#1}}
\newcommand{\ol}[1]{\overline{#1}}

\DeclareMathOperator{\End}{End}

\DeclareMathOperator{\Forg}{Forg}

\DeclareMathOperator{\Hom}{Hom}
\DeclareMathOperator{\Id}{Id}

\DeclareMathOperator{\Mod}{Mod}

\DeclareMathOperator{\Pro}{Pro}

\DeclareMathOperator{\add}{add}

\DeclareMathOperator{\biMod}{biMod}
\DeclareMathOperator{\codim}{codim}

\DeclareMathOperator{\coev}{coev}
\DeclareMathOperator{\colim}{colim}
\DeclareMathOperator{\comod}{comod}
\DeclareMathOperator{\bcomod}{\textbf{comod}}

\DeclareMathOperator{\ev}{ev}

\DeclareMathOperator{\im}{im}
\DeclareMathOperator{\imm}{Im}

\DeclareMathOperator{\modd}{mod}

\DeclareMathOperator{\rad}{rad}

\DeclareMathOperator{\zig}{zig}

\EnableBpAbbreviations
\author{James Macpherson}
\newtheorem{thm}{Theorem}[section]
\newaliascnt{lem}{thm}
\newtheorem{lem}[lem]{Lemma}
\aliascntresetthe{lem}
\newaliascnt{psn}{thm}
\newtheorem{psn}[psn]{Proposition}
\aliascntresetthe{psn}
\newaliascnt{cor}{thm}
\newtheorem{cor}[cor]{Corollary}
\aliascntresetthe{cor}
\newaliascnt{que}{thm}

\aliascntresetthe{que}
\newaliascnt{conj}{thm}

\newtheorem{thm*}{Theorem}
\aliascntresetthe{conj}
\theoremstyle{definition}
\newaliascnt{ex}{thm}

\aliascntresetthe{ex}
\newaliascnt{exs}{thm}

\aliascntresetthe{ex}
\newaliascnt{def}{thm}
\newtheorem{defn}[def]{Definition}
\aliascntresetthe{def}
\newaliascnt{not}{thm}
\newaliascnt{rmk}{thm}

\aliascntresetthe{rmk}
\newaliascnt{rmks}{thm}

\aliascntresetthe{rmks}
\usepackage{stmaryrd}
\title[2-Representations of Wide Finitary 2-Categories]{2-Representations and Associated Coalgebra 1-Morphisms for Locally Wide Finitary 2-Categories}

\begin{document}
	\bibliographystyle{alpha}
	
	\begin{abstract} We define locally wide finitary 2-categories by relaxing the definition of finitary 2-categories to allow infinitely many objects and isomorphism classes of 1-morphisms and infinite dimensional hom-spaces of 2-morphisms. After defining related concepts including transitive 2-representations in this setting, we provide a new method of constructing coalgebra 1-morphisms associated to transitive 2-representations of locally wide weakly fiat 2-categories, and demonstrate that any such transitive 2-representation is equivalent to a certain subcategory of the category of comodule 1-morphisms over the coalgebra 1-morphism. We finish the paper by examining two classes of examples of locally wide weakly fiat 2-categories: 2-categories associated to certain classes of infinite quivers, and singular Soergel bimodules associated to Coxeter groups with finitely many simple reflections.
		\end{abstract}
	
	\maketitle

\section{Introduction}

Finitary 2-representation theory is a powerful tool in representation theory, initially studied by Mazorchuk and Miemietz in the series of papers \cite{mazorchuk2011cell} through \cite{mazorchuk2016isotypic}, and then further developed by those authors and others in \cite{mackaay2016simple}, \cite{kildetoft2018simple} et al.. Its applications have included the representation theory of certain quotients of 2-Kac-Moody algebras (see \cite{mazorchuk2015transitive} and \cite{macpherson2020extension}) and Soergel bimodules corresponding to finite Coxeter groups (see \cite{kildetoft2018simple}).\par

However, the setup in which this work has been done has multiple finiteness restrictions imposed upon it. In particular, finitary 2-categories must have finitely many objects, and their hom-categories must each have finitely many isomorphism classes of indecomposable 1-morphisms and finite dimensional hom-spaces of 2-morphisms. In the earlier paper \cite{macpherson2020extension}, the author relaxed the finiteness restriction on the objects to study locally finitary 2-categories. This allowed new applications, including a wider class of quotients of 2-Kac-Moody algebras, but there are many interesting applications that remain outside that scope.\par

This paper seeks to relax all three restrictions, replacing each finite restriction with a countably infinite restriction, and imposing the requirement that the 2-category be locally Krull-Schmidt (which was guaranteed from the original definition of a finitary 2-category). These 2-categories, denoted \emph{locally wide finitary 2-categories}, open up many new potential applications for (wide) finitary 2-representation theory, including Soergel bimodules corresponding to any Coxeter group with a finite set of simple reflections.\par

The primary goal of this paper is to build up the necessary mathematical structure to prove \autoref{ProComodEquiv}: 
\begin{thm*} Any transitive 2-representation of a locally wide weakly fiat 1-category $\cC$ is equivalent to a certain subcategory of the comodule 1-morphism category over a coalgebra 1-morphism in a certain extension of $\cC$.\end{thm*} 

This allows any transitive 2-representation of a locally wide weakly fiat 2-category to be considered as an `internal' 2-representation, whose information is entirely encoded within the structure of $\cC$ itself. Once this result is established, we provide two classes of examples of locally wide finitary 2-categories, those of 2-categories associated to bound path algebras and of (singular) Soergel bimodules.\par

The structure of this paper is as follows. In \autoref{Pro2C} we recall various categorical and 2-categorical constructions we will be using in the paper, as well as constructing pro-bicategories, and proving in \autoref{StrPro2C} that the pro-bicategory of a strict 2-category is strict. \autoref{WF2Cats} defines locally wide finitary 2-categories and various accompanying concepts, such as wide finitary and transitive 2-representations. \autoref{AdelDef} defines the fan Adelman abelianisation of categories, 2-categories and 2-representations, which we will be using throughout the rest of the paper.\par

The main theorem of the paper lies in \autoref{CoalgDefn}, where for an object $S$ in a transitive 2-representation $\mbM$ of a locally wide weakly fiat 2-category $\cC$ we construct a corresponding coalgebra 1-morphism $[S,S]$ in the pro-2-category of the fan Adelman abelianisation of $\cC$. This leads to the main theorem, \autoref{ProComodEquiv}, where we construct a sub-category of the category of comodule 1-morphisms of $[S,S]$ which has the structure of a 2-representation of $\cC$, and is further equivalent to $\mbM$ as 2-representations of $\cC$.\par

\autoref{WFBPA} and \autoref{WFSoerg} concern themselves with examples of locally wide weakly fiat 2-categories that this theory applies to. The first example is that of 2-categories associated to pleasant classes of bound path algebras. In this case, we further prove in \autoref{basecoalgwidef} that the coalgebra 1-morphism $[S,S]$ is a 1-morphism in the locally wide weakly fiat 2-category $\cC$ itself, rather than simply living in the pro-2-category of the Adelman abelianisation. Further, in \autoref{stcellCAX} we classify all simple transitive 2-representations of these 2-categories as equivalent to cell 2-representations. The final section of the paper is devoted to proving  (in \autoref{SSoergWFiat}) that 2-categories of singular Soergel bimodules associated to Coxeter groups with finitely many simple reflections are locally wide fiat 2-categories, and hence amenable to the theory in the paper.

\newenvironment{acknowledgements} {\renewcommand\abstractname{Acknowledgements}\begin{abstract}} {\end{abstract}}

\begin{acknowledgements}
The author is very thankful to Vanessa Miemietz, who as PhD supervisor provided much support and advice during the PhD from which the material in this paper derives.\end{acknowledgements}

\tableofcontents

\section{(Pro-)2-Categories}\label{Pro2C}

\subsection{2-Categories}
We briefly recall the definitions for bicategories and 2-categories and give our notation. See \cite{johnson20212} Chapter 2 for a more detailed discussion of these concepts. Since our definition of wide finitary 2-categories given in \autoref{WFDefn} will generally constrain the sizes of any classes to countable sets, all classes in the following definitions are taken to be proper sets.
\begin{defn} A (small) \emph{bicategory} $\cC$ has the following data:\begin{itemize}
		\item A set of objects $\operatorname{Ob}(\cC)$, with elements generally denoted $\tti,\ttj,\dots\in\cC$.
		\item For each pair of object $\tti$ and $\ttj$, a category $\cC(\tti,\ttj)$ called the \emph{hom-category} from $\tti$ to $\ttj$. Objects of a hom-category are called \emph{1-morphisms} and generally denoted using uppercase Roman letters, commonly $X,Y,\dots$ or $F,G,\dots$. Morphisms of a hom-category are called \emph{2-morphisms} and are denoted using either lowercase Roman or lowercase Greek letters, with the choice generally made to increase clarity where possible. The morphism set between 1-morphisms $X$ and $Y$ is notated $\Hom_\cC(X,Y)$. We notate composition of 2-morphisms in a hom-category as $\beta\circ_V\alpha$. The identity 2-morphism of a 1-morphism $X$ is notated $\id_X$.
		\item For any three objects $\tti,\ttj$ and $\ttk$, a \emph{composition functor} $$\circ_{\tti,\ttj\ttk}:\cC(\ttj,\ttk)\times\cC(\tti,\ttj)\to\cC(\tti,\ttk).$$ For composable 1-morphisms $X$ and $Y$ and composable 2-morphisms $\alpha$ and $\beta$, we write $\circ_{\tti,\ttj,\ttk}(Y,X)=Y\circ X$ or $YX$ and $\circ_{\tti,\ttj,\ttk}(\beta,\alpha)=\beta\circ_H\alpha$.
		\item Distinguished \emph{identity 1-morphisms} $\bbon_\tti\in\cC(\tti,\tti)$.
		\item \emph{Associator} 2-isomorphisms $a_{X,Y,Z}:X\circ(Y\circ Z)\to(X\circ Y)\circ Z$ for any 1-morphisms $X,Y,Z$ where this makes sense, and left and right \emph{unital} 2-isomorphisms $\iota_X:\bbon_\tti\circ X\to X$ and $\rho_X: X\circ \bbon_\tti\to X$ respectively for any 1-morphism $X$ where this makes sense, all subject to certain axioms (see \cite{johnson20212} Section 2.1 for the full list of axioms).
	\end{itemize}
\end{defn}
\begin{defn} A bicategory $\cC$ is a \emph{2-category} if the associator and left and right unital 2-isomorphisms are all identity 2-morphisms.\end{defn}

\subsection{Pro-2-Categories}
We recall that the \emph{pro-category} $\Pro(\ceC)$ of a category $\ceC$ is a category whose objects are cofiltered diagrams of $\ceC$ (i.e. functors $I\to\ceC$ where $I$ is a cofiltered category). We denote such a diagram as $X=(X_i)_{i\in I}$ or $X=\underleftarrow{\lim}_{i\in I}X_i$ (dropping the labelling category where there is no confusion). Morphism sets are defined as $$\Hom_{\Pro(\ceC)}(X,Y)={\lim}_j\colim_i\Hom_{\ceC}(X_i,Y_j).$$\par

We can extend a functor $F:\ceC\to\ceD$ to a functor $\Pro(F):\Pro(\ceC)\to\Pro(\ceD)$. The action of $\Pro(F)$ on objects is given for $X=\underleftarrow{\lim}_i X_i$ by $\Pro(F)(X)=\underleftarrow{\lim}_i F(X_i)$. For morphisms, the map from $\Hom_{\Pro(\ceC)}(X,Y)$ to $\Hom_{\Pro(\ceD)}(FX,FY)$ is induced by the component maps $$F_{i,j}:\Hom_\ceC(X_i,Y_j)\to\Hom_\ceD(F(X_i),F(Y_j)).$$\par

We wish to construct a pro-bicategory $\Pro(\cC)$ of a bicategory $\cC$. Since we will be working with pro-2-categories of 2-categories for the main part of the paper, we need to take care and present our definitions explicitly. The objects of $\Pro(\cC)$ as the same as the objects of $\cC$, and we take the hom-categories in $\Pro(\cC)$ to be the pro-categories of the hom-categories of $\cC$; that is, $\Pro(\cC)(\tti,\ttj)=\Pro(\cC(\tti,\ttj))$ for any objects $\tti,\ttj\in\cC$. A 2-morphism therefore has the form $([\alpha_j])_j\in\prod\limits_{j\in J}((\coprod\limits_{i\in I}\Hom_\cC(X_i,Y_j))/\sim)$, where:\begin{itemize}
	\item For $\alpha_i:X_i\to Y_j$ and $\alpha_k:X_k\to Y_j$, $[\alpha_i]=[\alpha_k]$ if there exist 2-morphisms $\beta_i:X_l\to X_i$ and $\beta_k:X_l\to X_k$ in the diagram of $X$ with $\alpha_i\circ_V\beta_i=\alpha_k\circ_V\beta_k$. 
	\item For any 2-morphism $p:Y_j\to Y_k$ in the diagram of $Y$, $p[\alpha_j]=[\alpha_k]$.\end{itemize}
The vertical composition of two 2-morphisms is thus given by the standard composition in pro-categories.\par

For composition of 1-morphisms, given two composable 1-morphisms $X=\underleftarrow{\lim}_{i\in I}X_i$ and $Z=\underleftarrow{\lim}_{k\in K}Z_k$, we define the composite as $X\circ Z=\underleftarrow{\lim}_{i,k} X_i\circ Z_k$; that is, the diagram with 1-morphisms $X_i\circ Z_k$ and 2-morphisms $\alpha\circ_H\beta:X_i\circ Z_k\to X_j\circ Z_l$, where $\alpha:X_i\to X_j$ is a 2-morphism in the diagram of $X$ and $\beta:Z_k\to Z_l$ a 2-morphism in the diagram of $Z$. The identity 1-morphism of $\tti$ in $\Pro(\cC)$ is the `trivial' diagram consisting only of $\bbon_\tti$ and $\id_{\bbon_\tti}$.

It remains to define the horizontal composition of 2-morphisms. Given 1-morphisms $X=\underleftarrow{\lim}_{i\in I}X_i$, $Y=\underleftarrow{\lim}_{j\in J}Y_j$, $Z=\underleftarrow{\lim}_{k\in K}Z_k$ and $W=\underleftarrow{\lim}_{l\in L}W_l$, and 2-morphisms $([\alpha_k])_k:X\to Z$ and $([\beta_l])_l:Y\to W$, we can define the horizontal composition $([\alpha_k])\circ_H([\beta_l]):XY\to ZW$ to be $([\alpha_k])\circ_H([\beta_l])=([\alpha_k\circ_H\beta_l])$. It is straightforward to check that this is a well-defined element of $\Hom_{\Pro(\cC)}(XY,ZW)$. The identity 2-morphism of a 1-morphism $X=\underleftarrow{\lim}_{i\in I}X_i$ is $([\id_{X_i}])$.\par

Given these definitions, $\Pro(\cC)$ is a bicategory with the associativity 2-isomorphism $(X\circ Y)\circ Z\to X\circ (Y\circ Z)$ given by $([a_{X_i,Y_j,Z_k}])$, where the $a_{X_i,Y_j,Z_k}$ are the associativity 2-isomorphisms of $\cC$, and the right and left unital 2-isomorphisms given by $([\rho_{X_i}])_i:\bbon_\tti X\to X$ and $[(\iota_{X_i})]:X\bbon_\tti\to X$ respectively, where $\rho_{X_i}$ and $\iota_{X_i}$ are respectively the right and left unital 2-isomorphisms of $\cC$. It is straightforward but tedious to confirm that these obey the pentagon and unital axioms.\par

\begin{psn}\label{StrPro2C} In the setup of the prior paragraphs, if $\cC$ is a strict 2-category, then $\Pro(\cC)$ is a strict 2-category.
	\begin{proof} If $a_{X_i,Y_j,Z_k}=\id_{X_iY_jZ_k}$ for all $i$, $j$ and $k$, then it immediately follows that $([a_{X_i,Y_j,Z_k}])=([\id_{X_iY_jZ_k}])=\id_{XYZ}$ as required, and similarly for the left and right unital 2-morphisms.\end{proof}\end{psn}

When $\cC$ is a strict 2-category, we call $\Pro(\cC)$ a \emph{pro-2-category}.

\section{Wide Finitary 2-Categories}\label{WF2Cats}

Throughout this paper we take $\bbk$ to be an algebraically closed field of characteristic $0$.

\subsection{Initial Definitions}\label{WFDefn}
We now define the main 2-categories we will be studying in this paper. 

\begin{defn} A category $\ceC$ is \emph{wide finitary} if it is an additive $\bbk$-linear Krull-Schmidt category with countably many isomorphism classes of indecomposable objects and where the morphism sets are $\bbk$-vector spaces of countable dimension. We define the 2-category $\LA_\bbk^{wf}$ to have as objects wide finitary categories, as 1-morphisms $\bbk$-linear functors, and as 2-morphisms natural transformations.\end{defn}

\begin{defn} A 2-category $\cC$ is \emph{locally wide finitary} if: \begin{itemize}
		\item $\cC$ has countably many objects.
		\item For any objects $\tti,\ttj\in\cC$, $\cC(\tti,\ttj)\in\LA^{wf}_\bbk$.
		\item Horizontal composition is biadditive and $\bbk$-linear.
		\item For each object $\tti\in\cC$, the identity 1-morphism $\bbon_\tti$ is indecomposable.
	\end{itemize}
	
	If the 2-category only has finitely many objects, we refer to it just as a \emph{wide finitary} 2-category.\end{defn}

\begin{defn} Let $\cC$ be a locally wide finitary 2-category. If there exists a weak equivalence $-^*:\cC\to\cC^{\operatorname{co},\operatorname{op}}$ such that for any 1-morphism $X\in\cC(\tti,\ttj)$ there are natural 2-isomorphisms $\alpha:X\circ X^*\to\bbon_\tti$ and $\beta:\bbon_\ttj\to X^*\circ X$ with equalities $(\alpha\circ_H\id_X)\circ_V (\id_X\circ_H\beta)=\id_X$ and $(\id_{X^*}\circ_H\alpha)\circ_V(\beta\circ\id_{X^*})=\id_{X^*}$, then we say that $\cC$ is a \emph{locally wide weakly fiat 2-category}. We notate the inverse of $-^*$ by $\tensor*[^*]{-}{}$. If $-^*$ is an involution, then we say that $\cC$ is \emph{locally wide fiat}.\end{defn}

\begin{defn} Let $\cC$ be a locally wide finitary 2-category. A \emph{2-representation} of $\cC$ is a strict 2-functor from $\cC$ to $\cCat$. A \emph{wide finitary} 2-representation of $\cC$ is a strict 2-functor from $\cC$ to $\LA_\bbk^{\text{wf}}$. An \emph{abelian} 2-representation is a strict 2-functor from $\cC$ to the 2-category $\mathbf{AbCat}$ of abelian categories with additive functors and natural transformations.\end{defn}

We will make common use of the notation $\ceM=\coprod\limits_{\tti\in\cC}\mbM(\tti)$.

\subsection{Cells of Locally Wide Finitary 2-Categories}

\begin{defn} Let $\cC$ be a locally wide finitary 2-category. The set $S(\cC)$ of isomorphism classes of indecomposable 1-morphisms of $\cC$ has a multisemigroup structure with the operation given by $[F]*[G]=\{[H]\in S(\cC)|H\text{ is a direct summand of }FG\}$. We commonly abuse notation by writing $F$ in place of $[F]$.\end{defn}

We recall the Green's relations for multisemigroups, first defined in \cite{kudryavtseva2012multisemigroups} and stemming from the Green's relations for semigroups defined in \cite{green1951structure}. The left equivalence relation $\cL$ on $S(\cC)$ is given by $F\sim_\cL G$ if $S(\cC)*F=S(\cC)*G$; the right equivalence relation $\cR$ is defined analogously on the right, and the two-sided equivalence relation $\cJ$ is given by $F\sim_\cJ G$ if $$(F*S(\cC)*F)\cup (F*S(\cC))\cup (S(\cC)*F)=(G*S(\cC)*G)\cup (G*S(\cC))\cup (S(\cC)*G).$$ The equivalence classes of these relations are called $\cL$-, $\cR$- and $\cJ$-cells respectively.\par

There is an equivalent construction of these equivalence relations that we will utilise in this paper. We define partial orders $\leq_\cL$, $\leq_\cR$ and $\leq_\cJ$ on $S(\cC)$ by setting $F\leq_\cL G$ (respectively $F\leq_\cR G$, $F\leq_\cJ G$) if there exists some $H\in S(\cC)$ (respectively some $K\in S(\cC)$, some $M,N\in S(\cC)$) such that $G$ is a direct summand of $HF$ (respectively $G$ is a direct summand of $FK$, $G$ is a direct summand of $MFN$). The equivalence relation given by $F\sim_\cL G$ if and only if $F\leq_\cL G$ and $G\leq_\cL F$ is the same as that given above, and similarly for $\sim_\cR$ and $\sim_\cJ$.\par

We also mention here two further Green's relations: $\cH=\cL\cap\cR$ and $\cD$, which is the join of $\cL$ and $\cR$ in the poset of equivalence relations on $S(\cC)$. Working with cells in locally wide finitary 2-categories requires care, as we cannot assume that cells are finite in size or that there are only finitely many of them, even within a single hom-category. Indeed, we present a somewhat pathological example below:\par

Consider first a 2-category $\cD$ defined as follows: $\cD$ has only one object $*$. $\cD(*,*)$ is additive and $\bbk$-linear, with indecomposable 1-morphisms consisting of $\bbon_*$ and $F\llbracket z\rrbracket$ for $z\in\bbZ$, with composition defined by $F\llbracket z\rrbracket\circ F\llbracket y\rrbracket= F\llbracket z+y\rrbracket\oplus F\llbracket z+y+1\rrbracket$. For 2-morphisms we set $\Hom_\cD(F\llbracket z\rrbracket, F\llbracket y\rrbracket)\cong \delta_{yz}\bbk$ and extend additively. This clearly gives us a locally wide finitary 2-category. However, $F\circ F\llbracket z\rrbracket\cong F\llbracket z\rrbracket\oplus F\llbracket z+1\rrbracket$, while $F\llbracket -1\rrbracket \circ F\llbracket z+1\rrbracket\cong F\llbracket z\rrbracket\oplus F\llbracket z+1\rrbracket$. Therefore, $F\llbracket z\rrbracket \sim_L F\llbracket z+1\rrbracket$ for any $z$, and thus we only have two $\cL$-cells, one containing the identity 1-morphism and an infinitely large one containing all the $F\llbracket z\rrbracket$. We get similar results for $\cR$- and $\cJ$-cells.\par

Now construct the category $\cB$ as a disjoint union of countably many copies of $\cD$. It is clearly a locally wide finitary 2-category; however, by construction, it now has infinitely many infinitely large cells.\par

We might hope that imposing additional structure on these cells would provide nicer bounds for their cardinality. The natural structure to consider is that of strongly regular $\cJ$-cells, first defined for finitary 2-categories in \cite{mazorchuk2011cell}. It is certainly true that the proof \cite{mazorchuk2011cell} Proposition 28 b) is powerful enough to generalise to the locally wide finitary setting, giving the following result and allowing us to define strongly regular $\cJ$-cells and strongly regular locally wide finitary 2-categories in a similar fashion to \cite{mazorchuk2011cell}.

\begin{thm} Let $\cC$ be a locally wide finitary 2-category and let $\ceJ$ be a $\cJ$-cell of $\cC$ such that every $\cH$-cell of $\ceJ$ is non-empty. Let $\cL_\ceJ$, $\cR_\ceJ$, $\cD_\ceJ$ and $\cJ_\ceJ$ denote the restrictions of Green's relations to $\ceJ$. Then $\cL_\ceJ\circ\cR_\ceJ=\cR_\ceJ\circ\cL_\ceJ=\cD_\ceJ=\cJ_\ceJ$.\end{thm}

\begin{defn} A $\cJ$-cell $\ceJ$ of a locally wide finitary 2-category $\cC$ is \emph{strongly regular} if the $\cL$-cells of $\ceJ$ are incomparable under $\leq_\cL$ and if each $\cH$-cell of $\ceJ$ contains precisely one isomorphism class of indecomposables. The 2-category $\cC$ is \emph{strongly regular} if every $\cJ$-cell of $\cC$ is strongly regular.\end{defn}

Unfortunately, beyond the obvious restriction on the size of $\cH$-cells, even being strongly regular does not induce any further size limits on $\cL$-, $\cR$- or $\cJ$-cells. Consider a 2-category $\cD$ with a single object, and whose indecomposable 1-morphisms consist of the identity 1-morphism $\bbon$, and of 1-morphisms $F_{ij}$ for $i,j\in\bbZ$, where composition is defined by $F_{ij}\circ F_{kl}=F_{il}$. We have a $\cJ$-cell containing only the identity 1-morphism, which is trivially strongly regular, and then a $\cJ$-cell containing all the $F_{ij}$. Its $\cL$-cells are of the form $\ceL_i=\{F_{ix}|x\in\bbZ\}$, and its $\cR$-cells are of the form $\ceR_j=\{F_{yj}|y\in\bbZ\}$. The $\cL$-cells are clearly incomparable under the left order, and $\ceL_i\cap\ceR_j=\{F_{ij}\}$, and therefore the $\cJ$-cell in indeed strongly regular. However, the $\cJ$-cell and all its component $\cL$-and $\cR$-cells are infinite in size. We can again take $\cB$ to be countably many disjoint copies of $\cD$, to form a strongly regular locally wide finitary 2-category with infinitely many infinitely large $\cJ$-cells. For those interested in semi-group theory, $\cD$ is the 2-category induced by the rectangular band on $\bbZ\times\bbZ$ (see \cite{clifford1954bands}).\par

We also recall the standard definitions of a 2-ideal of a 2-category and an ideal of a 2-representation. 

\begin{defn} Given a $\bbk$-linear 2-category $\cC$, a \emph{left 2-ideal} $\cI$ of $\cC$ has the same objects as $\cC$, and for each pair $\tti,\ttj$ of objects an ideal $\cI(\tti,\ttj)$ of the 1-category $\cC(\tti,\ttj)$ stable under left horizontal multiplication with 1- and 2-morphisms of $\cC$. We similarly define \emph{right 2-ideals} and \emph{two-sided 2-ideals}. We call the latter simply \emph{2-ideals}.\end{defn}

\begin{defn} Given a $\bbk$-linear 2-category $\cC$ with a 2-ideal $\cI$ of $\cC$, the \emph{quotient 2-category} $\cC/\cI$ is defined as follows: the objects and 1-morphisms of $\cC/\cI$ are the same as those of $\cC$ and the hom-spaces between 1-morphisms are defined as $$\Hom_{\cC/\cI(\tti,\ttj)}(F,G)=\Hom_{\cC(\tti,\ttj)}(F,G)/\Hom_{\cI(\tti,\ttj)}(F,G).$$\end{defn}

\begin{defn} Given a locally wide finitary 2-category $\cC$ and a 2-representation $\mbM$ of $\cC$, an \emph{ideal} $\cI$ of $\mbM$ is a collection of ideals $\cI(\tti)\subseteq \mbM(\tti)$ which is closed under the action of $\cC$ in that for any morphism $f\in\cI(\tti)$ and any 1-morphism $F\in\cC$, $\mbM(F)(f)$ is a morphism in $\cI$ if it is defined.\end{defn}

\begin{defn} Given a locally wide finitary 2-category $\cC$, a 2-representation $\mbM$ of $\cC$ and an ideal $\cI$ of $\mbM$, we define the \emph{quotient 2-representation} $\mbM/\cI$ to be a 2-representation where for each object $\tti\in\cC$, $\mbM(\tti)$ and $(\mbM/\cI)(\tti)$ have the same objects, and where $\Hom_{\mbM/\cI}(X,Y)=\Hom_\mbM(X,Y)/(\cI(\tti)\cap\Hom_\mbM(X,Y))$ for any objects $X,Y\in(\mbM/\cI)(\tti)$.\end{defn}

\subsection{(Simple) Transitive and Cell 2-Representations}

For the duration of this section, we take $\cC$ to be a locally wide finitary 2-category.

\begin{defn} Let $\mbM$ be a wide finitary 2-representation of $\cC$. For a collection of objects $\Gamma\subseteq\operatorname{Ob}(\ceM)$, we define $\mbG_\mbM(\Gamma)$ to be the full sub-2-representation of $\mbM$ whose objects are the additive closure of objects of $\ceM$ of the form $\mbM(F)(X)$, where $X\in\Gamma\cap\mbM(\tti)$ for some $\tti\in\cC$ and $F\in\cC(\tti,\ttj)$ for some $\ttj\in\cC$. When $\Gamma=\{M\}$ is a single object of $\ceM$, we write $\mbG_\mbM(M)=\mbG_\mbM(\{M\})$.\end{defn}

\begin{defn} A wide finitary 2-representation $\mbM$ of $\cC$ is said to be \emph{transitive} if for any $M\in\ceM$, the 2-representation induced by $\mbG_\mbM(M)$ is equivalent to $\mbM$. Equivalently, $\mbM$ is transitive if for any $M,N\in\ceM$, $N$ is isomorphic to a direct summand of $\mbM(F)(M)$ for some 1-morphism $F$ of $\cC$. \end{defn}

We present a generalisation of \cite{mazorchuk2015transitive} Lemma 4.

\begin{lem}\label{WFMxId} Let $\cC$ be a locally wide finitary 2-category and let $\mbM$ be a transitive finitary 2-representation. There exists a unique maximal ideal $\cI$ of $\mbM$ which does not contain the identity morphism of any non-zero object.
	\begin{proof} We adapt the proof given in \cite{mazorchuk2015transitive} for our situation. The $\mbM(\tti)$ are additive categories and we can form the coproduct $\ceM=\coprod\limits_{\tti\in\cC}\mbM(\tti)$ in $\cCat$. Further, the structure of ideals of $\mbM$ does not depend on the number of indecomposable objects in any $\mbM(\tti)$, and so given an ideal $\mbK$ of $\mbM$ that does not contain any identity morphisms of non-zero objects, the coproduct $\cK=\coprod\limits_{\tti\in\cC} \mbK(\tti)$ is still an ideal of $\ceM$.\par
		
		Further, the $\mbM(\tti)$ are Krull-Schmidt, and thus by definition for $X$ an indecomposable object, $\End_\mbM(X)$ is a local algebra, and $\cK\cap\End_\mbM(X)$ is a proper ideal of $\End_\mbM(X)$. As the argument given in \cite{mazorchuk2015transitive} is a `pointwise' argument that considers the endomorphism ring of each indecomposable $X$ separately, and since an infinite sum of ideals is still an ideal, the proof generalises.
		
\end{proof}\end{lem}

\begin{defn} Let $\mbM$ be a transitive 2-representation of $\cC$, and let $\cI$ denote the maximal ideal as given in \autoref{WFMxId}. If $\cI=0$, then we say that $\mbM$ is a \emph{simple transitive} 2-representation. Given a transitive 2-representation $\mbM$ and its ideal $\cI$, we can form the simple transitive quotient 2-representation $\mbM^{{S}}=\mbM/\cI$, called the \emph{simple transitive quotient} of $\mbM$.\end{defn}

We define \emph{cell 2-representations} of $\cC$ analogously to previous papers (such as \cite{mazorchuk2011cell} for the finitary case and \cite{macpherson2020extension} for the locally finitary case). Let $\ceJ$ be a $\cJ$-cell of $\cC$ and let $\ceL\subseteq\ceJ$ be a $\cL$-cell. It is straightforward to show that there exists some object $\tti=\tti_{\ceL}$ of $\cC$ which is the domain of every $F\in\ceL$. For each $\ttj\in \cC$, set $$\mbN_\ceL(\ttj)=\add\{FX|F\in\cC(-,\ttj), X\in\ceL\}.$$ $\mbN_\ceL$ thus defines a map from the object set of $\cC$ into $\LA^{wf}_\bbk$. To make this a wide finitary 2-representation of $\cC$, we set $\mbN_\ceL(F)$ for a 1-morphism $F$ to be the functor defined by left composition with $F$, and $\mbN_\ceL(\alpha)$ for a 2-morphism $\alpha$ to be the natural transformation defined by left horizontal composition with $\alpha$.\par

We can generalise \cite{mazorchuk2015transitive} Lemma 3 in a similar fashion as the proof of \autoref{WFMxId} did for \cite{mazorchuk2015transitive} Lemma 4, and hence derive that $\mbN_\ceL$ has a unique maximal ideal $\cI_\ceL$ that does not contain $\id_F$ for any $F\in\ceL$.

\begin{defn} We define the \emph{cell 2-representation} $\mbC_\ceL$ as the simple transitive 2-representation given by $\mbN_\ceL/\cI_\ceL$.\end{defn}

Given the pathological examples given above, in general the cell 2-representations are indeed wide finitary 2-representations - the $\cL$-cells, regardless of any other comparable 1-morphisms, may contain infinitely many isomorphism classes of indecomposable 1-morphisms. As before, we notate the cell 2-representation of $\cC$ corresponding to a $\cL$-cell $\ceL$ by $\mbC_\ceL$.\par

\section{Adelman Abelianisation}\label{AdelDef}

\subsection{1-Categorical Construction}

Previous papers on (locally) finitary 2-representation utilised a version of abelianisation first defined by Freyd in \cite{freyd1966representations} and strictified in \cite{mackaay2016simple}. However, this method only constructs a legitimately abelian category when the original category has weak kernels or weak cokernels (depending on whether projective or injective abelianisation is performed). We present below a more powerful version of abelianisation. The basic construction is due to \cite{adelman1973abelian}, and we will present an equivalent but more complicated version to ensure we retain a 2-category once finished.

\begin{defn} Let $\ceC$ be an additive category. We construct the \emph{Adelman abelianisation} $\tilde{\ceC}$ as follows:\begin{itemize}
		\item The objects of $\tilde{\ceC}$ are quintuples $(Y,X,Z,\alpha,\beta)$, where $X,Y,Z$ are objects of $\ceC$ and $\alpha:Y\to X$ and $\beta:X\to Z$ are morphisms of $\ceC$.
		\item Morphisms of $\tilde{\ceC}$ are equivalence classes of triples $$(s,r,t):(Y,X,Z,\alpha,\beta)\to(Y',X',Z',\alpha',\beta')$$ where $s:Y\to Y'$, $r:X\to X'$ and $t:Z\to Z'$ are morphisms of $\ceC$, modulo those triples $(s,r,t)$ that satisfy a homotopy relation, explicitly triples $(s,r,t)$ for which there exist morphisms $p:X\to Y'$ and $q:Z\to X'$ such that $\alpha'p+q\beta=r$.
		\item Composition of triples is given by $(s,r,t)\circ(s',r',t')=(ss',rr',tt')$.
		\item Identity morphisms are of the form $(\id_Y,\id_X,\id_Z)$.\end{itemize}\end{defn}

\begin{defn} Let $\ceC$ be an additive category. We construct the \emph{fan Adelman abelianisation} $\widehat{\ceC}$ as follows:\begin{itemize}
		\item The objects of $\widehat{\ceC}$ are given by equivalence classes of sextuples of the form $(Y_i,X,Z_j,\alpha_i,\beta_j,k)_{i,j\in\bbZ^+}$ with $Y_i,X,Z_j\in \ceC$, $\alpha_i:Y_i\to X$, $\beta_j:X\to Z_j$ morphisms in $\ceC$ and $k\in\bbZ^+_0$. We require that for all $i,j>k$, $Y_i=Z_j=0$. Two sextuples are equivalent if they only differ in the value of $k$.
		\item Morphisms in $\widehat{\ceC}$ from $(Y_i,X,Z_j,\alpha_i,\beta_j,k)$ to $(Y_i',X',Z_j',\alpha_i',\beta_j',k')$ are equivalence classes of triples $(s_{i j}, r, t_{mn})_{i,j,m,n\in\bbZ^+}$ where $r:X\to X'$, $s_{ij}:Y_i\to Y_j'$ and $t_{mn}:Z_m\to Z_n'$ are morphisms in $\ceC$, with the equivalence relation being spanned by triples that satisfy a homotopy relation, explicitly triples $(s_{ij}, r, t_{mn})$ such that there exist $p_i:X\to Y_i'$ and $q_j:Z_j\to X'$ such that $\sum\limits_{i} \alpha_i' p_i+\sum\limits_j q_j\beta_j=r$.
		\item Composition of triples is given by $$(s_{ij}',r',t_{mn}')\circ(s_{ij},r,t_{mn})=(\sum\limits_l s_{lj}'s_{il}, r'r, \sum\limits_z t_{zn}'t_{mz}).$$
		\item Identity morphisms are of the form $(\delta_{ij}\id_{Y_i}, \id_X, \delta_{mn}\id_{Z_m})_{i,j,m,n\in\bbZ^+}$.
\end{itemize}\end{defn}

Similarly to the fan Freyd abelianisation from \cite{mackaay2016simple}, this can be thought of as a variation of the traditional Adelman abelianisation with multiple objects at the left and right, as in the diagram: \newline 
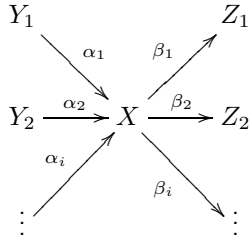
\xymatrix{ Y_1 \ar[rd]^{\alpha_1} & & Z_1 \\ Y_2 \ar[r]^{\alpha_2} & X \ar[ur]^{\beta_1} \ar[r]^{\beta_2} \ar[dr]_{\beta_i} & Z_2 \\ \vdots \ar[ur]^{\alpha_i} & & \vdots}  

While eventually both the $Y_i$ and the $Z_j$ will be zero, the minimal $i$ and $j$ where this occurs will in general not be identical. However, giving a single bound that is at least the larger of the two simplifies the already somewhat unwieldy notation.\par

In a similar fashion to the fan Freyd abelianisation found in \cite{mackaay2016simple}, $\widehat{\ceC}$ is additive and is equivalent to the traditional Adelman abelianisation via the assignment $$(Y_i,X,Z_j,\alpha_i,\beta_j,k)\mapsto \oplus Y_i \overset{\oplus \alpha_i}{\rightarrow} X \overset{\oplus\beta_j}{\rightarrow}\oplus Z_j$$ and is hence abelian. We can embed $\ceC$ into $\widehat{\ceC}$ via the assignment $X\mapsto (0,X,0,0,0,0)$ and $f:X\to Y\mapsto (0,f,0)$. As mentioned in \cite{adelman1973abelian}, its image is the full subcategory of injective-projectives of $\widehat{\cC}$. Let $I_\ceC:\ceC\to\widehat{\ceC}$ denote this canonical embedding of $\ceC$ into $\widehat{\ceC}$.\par

An important consequence of this construction is the following theorem due to \cite{adelman1973abelian}, using the aforementioned equivalence of categories between $\widehat{\ceC}$ and the original Adelman abelianisation:

\begin{thm}[\cite{adelman1973abelian} Theorem 1.14]\label{AbUnivExt} Let $\ceC$ be an additive 2-category, and let $F:\ceC\to\ceA$ be a additive functor into an abelian category $\ceA$. Then there is a unique (up to natural equivalence) exact functor $F^\diamond:\widehat{\ceC}\to \ceA$ such that $F^\diamond I_\ceC\cong F$.\end{thm}

If $\ceC$ has finite dimensional hom-spaces, then the dimension of $$\Hom_{\widehat{\ceC}}((Y_i,X,Z_j,\alpha_i,\beta_j,k),(Y_i',X',Z_j',\alpha_i',\beta_j',k'))$$ is bounded above by $$\prod\limits_{i\leq k, m\leq k'} \dim\Hom_\ceC(Y_i,Y_m')\cdot\dim\Hom_\ceC(X,X)\cdot\prod\limits_{j\leq k,n\leq k'}\dim\Hom_\ceC(Z_j,Z_n'),$$ and thus $\widehat{\ceC}$ also has finite dimensional hom-spaces. More generally, if $\ceC$ is a small additive category, then so is $\widehat{\ceC}$.

\subsection{2-Categorical Construction}

Let $\cC$ be a locally wide finitary 2-category. We define the \emph{fan Adelman abelianisation} $\widehat{\cC}$ of $\cC$ as follows:\begin{itemize}
	\item The objects of $\widehat{\cC}$ are the same as those of $\cC$.
	\item For any objects $\tti,\ttj\in\cC$, $\widehat{\cC}(\tti,\ttj)=\widehat{\cC(\tti,\ttj)}$, i.e. the hom-categories of $\widehat{\cC}$ are the fan Adelman abelianisations of the hom-categories of $\cC$.
	\item Composition of 1-morphisms is defined by $$(Y_i,X,Z_j,\alpha_i,\beta_j,k)\circ(Y_i',X',Z_j',\alpha_i',\beta_j',k')=(V_i,XX',W_j,\gamma_i,\delta_j,k+k'),$$ where:
	
	$$V_i=\begin{cases}Y_i\circ X', & i=1,\dots,k\\
	X\circ Y_{i-k}', & i=k+1,\dots,k+k' \\
	0, & \text{else,}\end{cases}$$ 
	
	$$W_j=\begin{cases}X\circ Z'_j, & j=1,\dots,k'\\
	Z_{j-k'}\circ X', & j=k'+1,\dots,k'+k \\
	0, & \text{else,}\end{cases}$$
	
	$$\gamma_i=\begin{cases} \alpha_i\circ_H\id_{X'}, & i=1,\dots,k\\
	\id_X\circ_H\alpha_{i-k}', & i=k+1,\dots,k+k'\\
	0, & \text{else,}\end{cases}$$
	
	$$\delta_j=\begin{cases} \id_X\circ_H\beta'_j, & j=1,\dots,k'\\
	\beta_{j-k'}\circ_H\id_{X'}, & j=k'+1,\dots,k'+k\\
	0, & \text{else.}\end{cases}$$
	
	\item Identity 1-morphisms are $(0,\bbon_\tti,0,0,0,0)$.
	\item Horizontal composition of 2-morphisms is defined component-wise.
\end{itemize}

The embedding of each 1-category $\ceC(\tti,\ttj)$ as the projective-injectives of $\widehat{\ceC}(\tti,\ttj)$ leads to the embedding of $\cC$ as a sub-2-category of $\widehat{\cC}$.\par

If $\cC$ is a locally wide finitary 2-category, let $\mbM$ be a wide finitary 2-representation of $\cC$. We define the \emph{fan Adelman abelianisation} $\widehat{\mbM}$ of $\mbM$ by setting $\widehat{\mbM}(\tti)=\widehat{\mbM(\tti)}$ for each object $\tti$ of $\cC$. This has the natural structure of a 2-representation of $\cC$ by component-wise action.\par

In addition, similar to \cite{mackaay2016simple}, we can make $\widehat{\mbM}$ a $\widehat{\cC}$ 2-representation by setting $\widehat{\mbM}((Y_i,X,Z_j,\alpha_i,\beta_j,k))(N_i,M,P_j,f_i,g_j,k')=(S_i,\mbM(X)M,T_j,h_i,l_j,k+k')$, where:

$$S_i=\begin{cases}\mbM(Y_i)M, & i=1,\dots,k\\
\mbM(X)N_{i-k}, & i=k+1,\dots,k+k' \\
0, & \text{else,}\end{cases}$$

$$T_j=\begin{cases}\mbM(X)P_j, & j=1,\dots,k'\\
\mbM(Z_{j-k'})M, & j=k'+1,\dots,k'+k \\
0, & \text{else,}\end{cases}$$

$$h_i=\begin{cases} \mbM(\alpha_i)_M, & i=1,\dots,k\\
\mbM(X)f_i, & i=k+1,\dots,k+k'\\
0, & \text{else,}\end{cases}$$

$$l_j=\begin{cases} \mbM(X)g_j, & j=1,\dots,k'\\
\mbM(\beta_{j-k'})_M, & j=k'+1,\dots,k'+k\\
0, & \text{else.}\end{cases}$$

\subsection{Beligiannis Abelianisation}

While we will mostly be using the fan Adelman abelianisation due to it being (relatively) simple to explicitly construct, there is another `universal' abelianisation that we will occasionally be referring to due to the process of its construction, originally defined in \cite{beligiannis2000freyd}.

\begin{defn} Given an additive category $\ceC$, we define the \emph{(fan) Beligiannis abelianisation} $\ul{\ol{\ceC}}=\ul{(\ol{\ceC})}(=\ol{(\ul{\ceC})})$ to be formed by taking the (fan) injective abelianisation of the (fan) projective abelianisation of $\ceC$.\end{defn}

By \cite{beligiannis2000freyd} Theorem 6.1 (1), this is equivalent to taking the (fan) projective abelianisation of the (fan) injective abelianisation of $\ceC$. Further, by \cite{beligiannis2000freyd} Theorem 6.1 (4) this has the same universal property as $\widehat{\ceC}$ does in \autoref{AbUnivExt}, and hence $\widehat{\ceC}$ and $\ul{\ol{\ceC}}$ are canonically equivalent, and since either fan Freyd abelianisation produces a strict 2-category when applied to a 2-category, this extends to a canonical equivalence of 2-categories between the abelianisations.

\section{Constructing the Coalgebra 1-Morphism}\label{CoalgDefn}

While we have a method of producing an abelianisation of a (locally) wide finitary (2-)category, we cannot immediately generalise the construction of the coalgebra 1-morphism from \cite{mackaay2016simple}, or indeed the algebra morphisms from \cite{etingof2016tensor} Section 7.8, as these both assume finiteness conditions that we lack. We thus need to expand to a larger setting.\par

Let $\cC$ be a locally wide finitary 2-category and let $\mbM$ be a transitive 2-representation of $\cC$. Choose $S\in\mbM(\tti)$. We define $\ceC_\tti=\coprod\limits_{\ttj\in\cC}\cC(\tti,\ttj)$. We also define the notation $\ceC_\tti(\ttj)=\cC(\tti,\ttj)$. We define a functor $\ev_S:\ceC_\tti\to\ceM$ by $\ev_S(F)=FS$ for a 1-morphism $F$ and for a 2-morphism $\alpha:F\to G$, $\ev_S(\alpha)=\alpha_S$. We let $\ev_{S,\ttj}$ denote the restriction and corestriction of $\ev_S$ to $\ceC_\tti(\ttj)$ and $\mbM(\ttj)$ respectively. By composing with the natural injection of $\ceM$ into $\widehat{\ceM}$, we can consider $\ev_{S,\ttj}$ to be a functor from $\ceC_\tti(\ttj)$ to $\widehat{\mbM(\ttj)}$. Since $\widehat{\mbM(\ttj)}$ is an abelian category, we can use the universal property of the Adelman abelianisation to extend $\ev_{S,\ttj}$ to an exact functor $\widehat{\ev_{S,\ttj}}:\widehat{\ceC_\tti}(\ttj)\to\widehat{\mbM(\ttj)}$. These then combine to give us a functor $\widehat{\ev_S}:\widehat{\ceC_\tti}\to\widehat{\ceM}$.

\begin{psn}\label{widehatEvnN} We can take $\widehat{\ev_S}$ to be evaluation at $S$.
	\begin{proof} We will show that evaluation at $S$ is an exact functor from $\widehat{\cC}(\tti,\ttj)$ to $\widehat{\mbM}(\ttj)$ for any $\ttj$. Then since $\widehat{\ev_{S,\ttj}}$ is the unique up to equivalence exact extension of $\ev_{S,\ttj}$, it must be equivalent to evaluation at $S$. For simplicity of notation we work in the non-fan Adelman abelianisation, since we are in essence dealing with a pair of 1-categories, and thus it is equivalent to the fan case.\par
		
		Consider some short exact sequence in $\widehat{\ceC}_\tti(\ttj)$,\newline \centerline{\xymatrix{ 0 \ar[d] & 0 \ar[d] & 0, \ar[d] \\ X_2 \ar[r]^{x_1} \ar[d]_{f_2} & X_1 \ar[r]^{x_2} \ar[d]_{f_1} & X_3 \ar[d]^{f_3} \\ Y_2 \ar[r]^{y_1} \ar[d]_{g_2} & Y_1 \ar[r]^{y_2} \ar[d]_{g_1} & Y_3 \ar[d]^{g_3} \\ Z_2 \ar[r]^{z_1} \ar[d] & Z_1 \ar[r]^{z_2} \ar[d] & Z_3 \ar[d] \\ 0 & 0 & 0}} which we also notate as $0 \to X \overset{f}{\to} Y \overset{g}{\to} Z \to 0$ for brevity. Since $\widehat{\ceC}_\tti(\ttj)$ is an abelian category and $f$ is monic, $\imm f=Y$ and so $f=\ker g$. But by \cite{adelman1973abelian} Theorem 1.1, we thus have an explicit construction for $\ker g$ (up to isomorphism), namely\newline \centerline{\xymatrix{Y_2\oplus Z_2 \ar[r]^{\begin{pmatrix} y_1 & 0 \\ g_2 & -1\end{pmatrix}} \ar[d]_{(1,0)} & Y_1\oplus Z_2 \ar[r]^{\begin{pmatrix} g_1 & -z_1 \\ y_2 & 0\end{pmatrix}} \ar[d]_{(1,0)} & Z_1\oplus Y_3 \ar[d]^{(0,1)} \\ Y_2 \ar[r]^{y_1} & Y_1 \ar[r]^{y_2} & Y_3.}}
		
		But the evaluation of this at $S$ is\newline \centerline{\xymatrix{Y_2S\oplus Z_2S \ar[rr]^{\begin{pmatrix} (y_1)_S & 0 \\ (g_2)_S & (-1)_S\end{pmatrix}} \ar[d]_{(1,0)} & & Y_1S\oplus Z_2S \ar[rr]^{\begin{pmatrix} (g_1)_S & (-z_1)_S \\ (y_2)_S & 0\end{pmatrix}} \ar[d]_{(1,0)} & & Z_1S\oplus Y_3S \ar[d]^{(0,1)} \\ Y_2S \ar[rr]^{(y_1)_S} & & Y_1S \ar[rr]^{(y_2)_S} & & Y_3S}} and since $\widehat{\mbM}(\ttj)$ is an Adelman abelianisation, this is the kernel of $g_S$. Thus $f_S$ is monic and since $\widehat{\mbM}(\ttj)$ is abelian, $\im f_S=\ker g_S$.\par
		
		It remains to show that $g_S$ is epic. But $g$ is epic, and since $\widehat{\ceC}_\tti(\ttj)$ is abelian, it is thus a cokernel of some morphism $h:W\to Y$. But again using \cite{adelman1973abelian} Theorem 1.1 and a similar procedure to above, we derive that $g_S$ is the cokernel of $h_S:WS\to YS$, and thus is epic. Therefore $f_S$ is monic, $g_S$ is epic and $\im f_S=\ker g_S$, and hence the sequence $0\to XS\overset{f_S}{\to} YS \overset{g_S}{\to} ZS\to 0$ is exact in $\widehat{\mbM}(\ttj)$, and the result follows.
\end{proof}\end{psn}

Before continuing, we wish to define the action of the pro-2-category $\Pro(\widehat{\cC})$ on $\Pro(\widehat{\ceM})$. We define this component-wise: since we have the action of $\widehat{\cC}(\tti,\ttj)$ on $\widehat{\mbM}(\tti)$, i.e. a bifunctor $\ev_{-}(-):\widehat{\cC}(\tti,\ttj)\times\widehat{\mbM}(\tti)\to\widehat{\mbM}(\ttj)$, we can take the pro-functor $\Pro(\ev_{-}(-))$ as the action by using a similar process to the definition of a pro-2-category in \autoref{Pro2C}. In particular, keeping the same $S$ as above and taking $X=\underleftarrow{\lim}_i X_i\in\Pro(\ceC_\tti(\ttj))$, we have that $XS=\underleftarrow{\lim}_i X_iS\in\Pro(\widehat{\mbM}(\ttj))$.

\begin{psn}\label{ProEvnN} $\Pro(\widehat{\ev_S})$ is evaluation at $S$.
	\begin{proof} Let $X=\underleftarrow{\lim}_i X_i\in\Pro(\widehat{\ceC_\tti}(\ttj))$. Then $$\Pro(\widehat{\ev_{S,\ttj}})(X)=\underleftarrow{\lim}_i(\widehat{\ev_{S,\ttj}}(X_i))=\underleftarrow{\lim}_i X_iS= XS$$ by \autoref{widehatEvnN}, and the result follows.\end{proof}\end{psn}

We recall here a pair of results from \cite{grothendieck1972prefaisceaux}, though we give the dual versions thereof:

\begin{lem}[\cite{grothendieck1972prefaisceaux} Proposition 8.11.4]\label{GV8_11_4} Let $\ceC$ be a category equivalent to a small category and let $\ceD$ be a category. Then a functor $F:\ceC\to\ceD$ has a pro-adjoint if and only if it is right exact.\end{lem}

\begin{lem}[\cite{grothendieck1972prefaisceaux} Proposition 8.11.2]\label{GV8_11_2} Let $\ceC$ and $\ceD$ be categories. Then a functor $F:\ceC\to\ceD$ has a pro-adjoint if and only if $\Pro(F)$ has a left adjoint.\end{lem}

\begin{psn} $\Pro(\widehat{\ev_S})$ has a left adjoint, denoted $[S,-]:\Pro(\widehat{\ceM})\to\Pro(\widehat{\ceC_\tti})$.
	\begin{proof} By construction $\widehat{\ev_S}$ is an exact functor, which by \autoref{GV8_11_4} is equivalent to $\widehat{\ev_S}$ having a pro-adjoint, which by \autoref{GV8_11_2} is equivalent to $\Pro(\widehat{\ev_S})$ having a left adjoint, as required.\end{proof}\end{psn}

We are especially interested in $[S,S]$, since we will show that $[S,S]$ has the structure of a coalgebra 1-morphism in $\Pro(\widehat{\ceC_\tti}(\tti))$ when $\cC$ is sufficiently pleasant. An element $m$ of $\Hom_{\Pro(\widehat{\ceC_\tti})}(X,Y)$ acts as a function $\Hom_{\Pro(\widehat{\ceC_\tti})}([S,S],X)\to\Hom_{\Pro(\widehat{\ceC_\tti})}([S,S],Y)$ via composition. But using the adjunction isomorphism, this can also be considered as a function $\Hom_{\Pro(\widehat{\mbM}(\tti))}(S,XS)\to \Hom_{\Pro(\widehat{\mbM}(\tti))}(S,YS)$ via composition with $m_S$.

\begin{psn}\label{AdjnProC} If $\cC$ is a locally wide weakly fiat 2-category and $G,H\in\Pro(\widehat{\cC})$, $\Hom_{\Pro(\widehat{\cC})}(FG,H)\cong\Hom_{\Pro(\widehat{\cC})}(G,F^*H)$ for any $F\in\cC$.
	\begin{proof} For a 1-morphism $F\in\cC(\tti,\ttj)$, the evaluation and coevaluation 2-morphisms $\alpha:FF^*\to\bbon_\ttj$ and $\beta:\bbon_\tti\to F^*F$ still exist in $\Pro(\widehat{\cC})$ and $\End_{\Pro(\widehat{\cC})}(G)=\End_\cC(G)$ for any $G\in\cC$, thus it remains the case that $(\alpha\circ_H\id_F)\circ_V(\id_F\circ_H\beta)=\id_F$ and $(\id_{F^*}\circ_H\alpha)\circ_V(\beta\circ_H\id_{F^*})=\id_{F^*}$, and thus $F$ and $F^*$ still form an internal adjunction as required.\end{proof}\end{psn}

An immediate consequence is the following:

\begin{cor}\label{AdjnProM} Let $\cC$ be a locally wide weakly fiat 2-category, $\mbM$ a wide finitary 2-representation of $\cC$ and $X,Y\in\Pro(\widehat{\ceM})$. Then $$\Hom_{\Pro(\widehat{\ceM})}(FX,Y)\cong\Hom_{\Pro(\widehat{\ceM})}(X,F^*Y)$$ for any $F\in\cC$\end{cor}

\begin{psn}\label{ProCoalgExist} $[S,S]$ has the structure of a coalgebra 1-morphism in $\Pro(\widehat{\ceC_\tti})$.
	\begin{proof} We construct comultiplication and counit 2-morphisms for $[S,S]$ analogously to the proof of \cite{mackaay2016simple} Lemma 4.3. For the comultiplication, by construction of the adjunction we have an isomorphism of hom-spaces \begin{align*}\Hom_{\Pro(\widehat{\ceC_\tti})}([S,S],[S,S])& \cong\Hom_{\Pro(\widehat{\mbM}(\tti))}(S,\Pro(\widehat{\ev_S})[S,S])\\ & =\Hom_{\Pro(\widehat{\mbM}(\tti))}(S,[S,S]S),\end{align*} with the equality following from \autoref{ProEvnN}. Let $\coev_S:S\to [S,S]S$ be the image of $\id_{[S,S]}$ under this isomorphism. We can thus form the composition $$([S,S]\coev_S)\circ\coev_S\in\Hom_{\Pro(\widehat{\mbM}(\tti))}(S,[S,S][S,S]S).$$ But again by the adjunction isomorphism, $$\Hom_{\Pro(\widehat{\mbM}(\tti))}(S,[S,S][S,S]S)\cong\Hom_{\Pro(\widehat{\ceC_\tti})}([S,S],[S,S][S,S]).$$ We take the coevaluation 2-morphism $\delta_S$ to be the image of $([S,S]\coev_S)\circ\coev_S$.\par
		
		For the counit 2-morphisms, we have the adjunction isomorphism $$\Hom_{\Pro(\widehat{\mbM}(\tti))}(S,S)=\Hom_{\Pro(\widehat{\mbM}(\tti))}(S,\bbon_\tti S)\cong\Hom_{\Pro(\widehat{\ceC_\tti})}([S,S],\bbon_\tti).$$ We thus take the counit 2-morphism $\epsilon_S$ to be the image of $\id_S$ under this isomorphism.\par
		
		Showing that $\delta_S$ and $\epsilon_S$ satisfy the coalgebra axioms is a collection of straightforward computations, giving the result.\end{proof}\end{psn}

\begin{psn} Let $T\in\mbM(\ttj)$. Then $[S,T]$ is a comodule 1-morphism over $[S,S]$.
	
	\begin{proof} The proof is analogous to that of \autoref{ProCoalgExist}.\end{proof} \end{psn}

We denote the category of comodule 1-morphisms over $[S,S]$ by $\comod_{\Pro(\widehat{\cC})}([S,S])$, which we abbreviate to $\comod([S,S])$ when it does not cause confusion. This can also be considered as a 2-representation of $\cC$ and $\Pro(\widehat{\cC})$ in a canonical fashion; when we do so, we denote the corresponding 2-representation by $\bcomod_{\Pro(\widehat{\cC})}([S,S])$ or $\bcomod([S,S])$. This allows us to define a functor $\Theta:\ceM\to\comod_{\Pro(\widehat{\cC})}([S,S])$ given on objects by $T\mapsto [S,T]$ and on morphisms by $f\mapsto [S,f]$. We denote by $\operatorname{Forg}_S:\comod_{\Pro(\widehat{\cC})}([S,S])\to \Pro(\widehat{\ceC})$ the canonical forgetful functor.

\begin{psn} $\Theta$ is indeed a functor.
	
	\begin{proof} It suffices to show that $[S,f]:[X,T]\to[S,T']$ is a morphism in $\comod([S,S])$ for any morphism $f:T\to T'$ in $\ceM$. Specifically, the diagram \newline \centerline{\xymatrix{ [S,T] \ar[rr]^{\rho_T} \ar[d]_{[S,f]} & & [S,T][S,S] \ar[d]^{[S,f]\circ_H\id_{[S,S]}} \\ [S,T'] \ar[rr]_{\rho_{T'}} & & [S,T'][S,S]}} needs to commute. We will show this by showing that the images of the sides under the adjunction isomorphism are equal.\par
		
		For notation, let $\coev_{S,T}$ denote the image of $\id_{[S,T]}$ under its adjunction isomorphism. Letting $\eta$ be the unit of the adjunction and $\sigma$ the counit, we have that $\eta_T=\coev_{S,T}$, and given $\alpha\in\Hom_{\Pro(\widehat{\ceC})}([S,T],F)$, the image of $\alpha$ under the adjunction isomorphism is given by $\alpha_S\circ\eta_T=\alpha_S\circ\coev_{S,T}$. Similarly, the image of $f\in\Hom_{\Pro(\widehat{\ceM})}(T,FS)$ is $\sigma_{[S,T]}\circ_V[S,f]$.\par
		
		Under the transferral of the action in the previous paragraph, $([S,f]\circ_H\id_{[S,S]})\circ_V\rho_T$ maps to $$([S,f]\circ_H\id_[S,S])_S\circ[S,T]\coev_S\circ\coev_{S,T}=[S,f]_{[S,S]S}\circ[S,T]\coev_S\circ\coev_{S,T}.$$ We wish to show that the diagram \newline \centerline{\xymatrix{T \ar[r]^{\coev_{S,T}} \ar[d]_f & [S,T]S \ar[rr]^{[S,T]\coev_S} \ar[d]^{[S,f]_S} & & [S,T][S,S]S \ar[d]^{[S,f]_{[S,S]S}} \\ T' \ar[r]_{\coev_{S,T'}} & [S,T']S \ar[rr]_{[S,T']\coev_S} & & [S,T'][S,S]S}} commutes. The right hand square does in fact commute, since $\mbM([S,f])$ is a natural transformation. For the left hand square, the image of $\coev_{S,T'}\circ f=\eta_{T'}\circ f$ under the adjunction isomorphism is $\sigma_{[S,T']}\circ_V[S,\eta_{T'}\circ f]$. But by the triangle identities for adjunctions, $\sigma_{[S,T']}\circ_V[S,\eta_{T'}\circ -]=[S,-]$, and therefore the image of $\coev_{S,T'}\circ f$ is $[S,f]$, and the reverse isomorphism takes $[S,f]$ to $[S,f]_S\circ\coev_{S,T}$, giving the required isomorphism.\par

		Thus $$[S,f]_{[S,S]S}\circ[S,T]\coev_S\circ\coev_{S,T}=[S,T']\coev_S\circ\coev_{S,T'}\circ f.$$ Similarly, $\rho_{T'}\circ_V[S,f]$ maps to $(\rho_{T'})_S\circ[S,f]_S\circ\coev_{S,T}$, which again by the commutativity of the left square above is equal to $(\rho_{T'})_S\circ\coev_{S,T'}\circ f$. We thus wish to show that $$(\rho_{T'})_S\circ\coev_{S,T'}\circ f=[S,T']\coev_S\circ\coev_{S,T'}\circ f.$$ But using the adjunction isomorphism and the definition of $\rho_{T'}$, we have that $$[S,T']\coev_S\circ\coev_{S,T'}\mapsto\rho_T'\mapsto(\rho_{T'})_S\circ\coev_{S,T'},$$ and the result follows.\par\end{proof}\end{psn}

\begin{defn} We define the \emph{image of $\ceM$ in $\comod([S,S])$} to be the full subcategory  $[S,\ceM]$ of $\comod([S,S])$ whose objects consist of those objects isomorphic to $\{([S,T],\rho_T)\in\comod([S,S])|T\in\ceM\}$. Similarly, we define $[S,\ceM(\tti)]$ to be the full subcategory of $\comod([S,S])$ whose objects consist of those objects isomorphic to $\{([S,T],\rho_T)\in\comod([S,S])|T\in\ceM(\tti)\}$.\end{defn}

By definition the image of $\Theta$ is contained in $[S,\ceM]$. Further, $[S,\ceM]$ has the structure of a locally wide finitary 2-representation $[S,\mbM]$ of $\cC$ by setting $[S,\mbM](\tti)=[S,\ceM(\tti)]$, $([S,\mbM](F))([S,T])=[S,\mbM(F)T]$ and $([S,\mbM](\alpha))_T=[S,\mbM(\alpha)_T]$. For the rest of this section, assume that $\cC$ is a locally wide weakly fiat 2-category.

\begin{psn} $\Theta$ defines a morphism of 2-representations of $\cC$ between $\mbM$ and $\bcomod([S,S])$.
	
	\begin{proof} We mirror the proof of \cite{mackaay2016simple} Lemma 4.4. We will start by showing that $[S,FT]\cong F[S,T]$ for any $F\in\cC(\tti,\ttj)$. Indeed, for any $G\in\Pro(\widehat{\cC})(\tti,\ttj)$ we have \begin{align*} \Hom_{\Pro(\widehat{\cC})}([S,FT],G) &\cong \Hom_{\Pro(\widehat{\ceM})}(FT,GS)\\ 
		& \cong\Hom_{\Pro(\widehat{\ceM})}(T,F^*GS)\\ 
		& \cong\Hom_{\Pro(\widehat{\cC})}([S,T],F^*G)\\
		& \cong\Hom_{\Pro(\widehat{\cC})}(F[S,T],G)\end{align*} using \autoref{AdjnProC} and \autoref{AdjnProM}, from which the claim follows.\par
		
		The rest of the proof is a direct generalisation of the remainder of the proof of \cite{mackaay2016simple} Lemma 4.4.\end{proof}\end{psn}

\begin{psn} For any $F\in\Pro(\widehat{\ceC_\tti})$, $[S,FS]\cong F[S,S]$ is a cofree $[S,S]$-coalgebra; that is, for any comodule $X\in\comod([S,S])$, $$\Hom_{\comod([S,S])}(X,F[S,S])\cong\Hom_{\Pro(\widehat{\cC})}(X,F).$$
	
	\begin{proof} We construct the relevant adjunction. We will again use the forgetful functor $\Forg:\comod([S,S])\to\Pro(\widehat{\ceC_\tti})$. Let $\ev_{[S,S]}:\Pro(\widehat{\ceC_\tti})\to\comod([S,S])$ be defined for a 1-morphism $F\in \Pro(\widehat{\ceC_\tti})$ as $\ev{[S,S]}(F)=(F[S,S], \id_F\circ_H\delta_S)$ and for a 2-morphism $\beta:F\to G$ as $\ev_{[S,S]}(\beta)=\beta\circ_H\id_{[S,S]}$. This is indeed a functor because the diagram\newline \centerline{\xymatrix{F[S,S] \ar[rr]^{\id_F\circ_H\delta_S} \ar[d]_{\beta\circ_H\id_{[S,S]}} & & F[S,S][S,S] \ar[d]^{\beta\circ_H\id_{[S,S][S,S]}} \\ G[S,S] \ar[rr]_{\id_G\circ_H \delta_S} & & G[S,S][S,S]}} clearly commutes. We claim that $\ev_{[S,S]}$ is a right adjoint of $\Forg$, from which the result follows immediately.\par
		
		To see this is indeed an adjunction, we will define the unit and counit. The unit $\eta:\Id_{\comod([S,S])}\to \ev_{[S,S]}\circ \Forg$ is defined by $\eta_{(X,\rho_X)}=\rho_X$, and the counit by $\sigma:\Forg\circ\ev_{[S,S]}\to\Id_{\Pro(\widehat{\ceC})}$ by $\sigma_F=\id_F\circ_H\epsilon_S$. The left triangle identity is thus expressed for $X\in\comod([S,S])$ by \xymatrix{X \ar[rr]^{\Forg(\rho_X)} & & X[S,S] \ar[rr]^{\id_X\circ_H\epsilon_S} & & X}, which is the composite $(\id_X\circ_H\epsilon_S)\circ_V\rho_X$, which is $\id_X$ by the comodule axioms. The right triangle identity is expressed for $F\in\Pro(\widehat{\ceC_\tti})$ by\newline \centerline{\xymatrix{F[S,S] \ar[rr]^{\id_F\circ_H\delta_S} & & F[S,S][S,S] \ar[rrr]^{\id_F\circ_H\epsilon_S\circ_H\id_{[S,S]}} & & & F[S,S]},} i.e. the composite $$(\id_F\circ_H(\epsilon_S\circ_H\id_{[S,S]}))\circ_V(\id_F\circ_H\delta_S)=\id_F\circ_H((\epsilon_S\circ_H\id_{[S,S]})\circ_V\delta_S),$$ which is equal to $\id_F\circ_H\id_{[S,S]}$ by the coalgebra axioms. We are thus done.\end{proof}\end{psn}

\begin{thm}\label{ProComodEquiv} $\Theta$ define an equivalence of 2-representations of $\cC$ between $\mbM$ and $[S,\mbM]$.
	
	\begin{proof} We mirror the proof of \cite{mackaay2016simple} Theorem 4.7.  By definition, $\Theta$ is essentially surjective when corestricted to $[S,\ceM]$, and it remains to show that it is fully faithful. To start, consider $FS,GS\in\ceM$ for $F,G\in\cC$. Then we have \begin{align*} \Hom_{\comod([S,S])}(\Theta(FS),\Theta(GS)) & \cong \Hom_{\comod([S,S])}(F[S,S],G[S,S]) \\ &\cong\Hom_{\Pro(\widehat{\cC})}(F[S,S],G) \\ &\cong\Hom_{\Pro(\widehat{\cC})}([S,S],F^*G) \\ &\cong \Hom_{\ceM}(S,F^*GS) \\ &\cong\Hom_\ceM(FS,GS).\end{align*} But now since $\mbM$ is a transitive 2-representation, for any $T_1,T_2\in\ceM$ there exist some 1-morphisms $F,G\in\cC$ such that $T_1$ is a direct summand of $FS$ and $T_2$ is a direct summand of $GS$. Thus by pre- and post-composing with injection and projection morphisms, and using that $\Theta$ preserves biproducts, we derive that $$\Hom_{\comod([S,S])}(\Theta(T_1),\Theta(T_2))\cong\Hom_\ceM(T_1,T_2).$$ Hence $\Theta$ is fully faithful, and the result follows.\end{proof}\end{thm}

\section{Application: Bound Path Algebras}\label{WFBPA}

We being by recalling the construction of bound path algebras. Let $\Gamma$ be a quiver, and let $\bbk\Gamma$ denote its path algebra. We let $\bbk\Gamma_i$ denote the ideal of $\bbk\Gamma$ generated by all paths of length $i$. A \emph{bound path algebra} is then a quotient $\bbk\Gamma/I$ of $\bbk\Gamma$ by some ideal $I$ such that there exists some integer $n\geq 2$ with $\bbk\Gamma_n\subseteq I\subseteq\bbk\Gamma_2$. To minimise the chance of confusion, we also remind the reader of the definition of a \emph{local} algebra (an algebra with a unique maximal left ideal which is equal to the Jacobson radical), and that this is unrelated to a \emph{locally [property]} 2-category (a 2-categories whose hom-categories are all [property] categories).\par

Let $A$ be a connected non-unital self-injective bound path algebra over $\bbk$. We notate by $V_A$ the set of vertices of the underlying quiver of $A$, and without loss of generality we assume that $V_A\subseteq\bbZ$. For our purposes, we also assume that every element of $V_A$ has finite total vertex degree. In this case, $A$ has the following property: $A$ has an orthogonal set of primitive idempotents $\{e_i|i\in V_A\}$ such that $A=\bigoplus\limits_{i,j\in V_A} e_iAe_j$, and such that $Ae_i$ and $e_jA$ are finite dimensional over $\bbk$ for all $i,j\in V_A$. Without loss of generality the $e_i$ are the paths of length 0.\par

As an example, consider the algebra $A_{\zig}$ stemming from the infinite zigzag quiver; that is, the quiver with vertices labelled by the integers and arrows $a_i:i\to i+1$ and $b_i:i+1\to i$ such that $a^2=b^2=0$ and $ab=ba$. In this case, $V_A=\bbZ$, and $A_{\zig}$ is generated by $e_i$, $a_i$ and $b_i$ subject to the following conditions:\begin{itemize}
	\item $\{e_i|i\in\bbZ\}$ is a set of idempotents fitting the above property;
	\item $a_i\in e_{i+1}A_{\zig}e_i$; $b_i\in e_iA_{\zig}e_{i+1}$ for all $i$;
	\item $a^2=b^2=0$ and $ab=ba$.\end{itemize}

There is a pleasant classification of the finite dimensional projective modules of $A$:

\begin{psn}\label{BPAProj} Every indecomposable projective module of $A$-$\modd$ is of the form $Ae_i$ for some $i\in A$.
	\begin{proof} It is clear that each $Ae_i$ is an indecomposable projective module, and thus it suffices to show that they exhaust the indecomposable projective modules. Let $M$ be some finite dimensional $A$ module. Since it is finite dimensional, it has a composition series $\{0\}=M_0\subset M_1\subset\dots\subset M_n=M$. Let $L_t=M_t/M_{t-1}$ be the $t$th (simple) composition factor. Let $i_t\in V_A$ be such that $e_{i_t}L_t\neq 0$. Then we have surjections $Ae_{i_t}\twoheadrightarrow L_t$ and $M_t\twoheadrightarrow L_t$. By the universal property of projective modules, this implies there is a homomorphism $Ae_{i_t}\to M_t$ making the resulting diagram commutative.\par
		
		These morphisms compile to form a surjection $q:\bigoplus\limits_{t=1}^n Ae_{i_t}\twoheadrightarrow M$. In particular, if $M$ is projective, then $q$ is a split epimorphism by the universal property of projective modules (applied to $M$), and hence $M$ is a direct summand of $\bigoplus\limits_{t=1}^n Ae_{i_t}$. The result follows.\end{proof}\end{psn}

In particular, \autoref{BPAProj} implies that every indecomposable projective module in $A$-$\modd$ is finitely generated and finite dimensional. This allows us to use standard module theoretic results to derive the following:

\begin{psn} Every simple module $S_i$ in $A$-$\modd$ is the quotient of some $Ae_i$ by its (unique) maximal submodule $R_i$.\end{psn}

\begin{psn} The maximal submodule $R_i$ of $Ae_i$ is generated by the equivalence classes of paths of length $1$. 
\begin{proof} Let $R_i'$ be the submodule of $Ae_i$ generated by the equivalence classes of paths of length $1$. By definition $R_i'$ is a submodule, and since $A$ is a bound path algebra, it is a proper submodule. Further, let $r\in Ae_i\setminus R_i'$. Then by construction $r= me_i+r'$, where $m\in\bbk$ and $r'$ does not have $e_i$ as a summand. Let $M$ be a submodule of $Ae_i$ such that $R_i'\cup\{r\}\subseteq M$. Then by construction $r'\in R_i'$, and hence $m^{-1}(r-r')=e_i\in M$, and hence $M=Ae_i$. The result follows.\end{proof}\end{psn}

This implies that $S_i=\bbk\hat{e_i}$, where the action of $A$ is given by $e_i(\hat{e_i})=\hat{e_i}$ and $p\hat{e_i}=0$ for any equivalence class of paths $p\neq e_i$ in $A$.\par

As an example, in $A_{\zig}$-$\modd$, the indecomposable projective $A_{\zig}e_i$ is four-dimensional as a $\bbk$-vector space, with a canonical basis $\{e_i,a_i,b_{i-1}, a_ib_i=b_{i-1}a_{i-1}\}$.\par

From the prior paragraphs, it follows that the indecomposable projective $(A$-$A)$-bimodules are of the form $Ae_i\otimes_\bbk e_jA$ for $i,j\in V_A$. For compactness of notation, we set $A_{ij}=Ae_i\otimes_\bbk e_jA$. This allows us to construct a 2-category $\cC_A$ as follows:\begin{itemize}
	
	\item $\cC_A$ has one object, which we associate with (a small category equivalent to) $A$-$\modd$.
	
	\item 1-morphisms are isomorphic to direct summands of direct sums of the identity and of functors isomorphic to tensoring with $Ae_i\otimes_\bbk e_jA$ for $i,j\in V_A$. For compact notation, we set $F_{ij}= Ae_i\otimes_\bbk e_jA\otimes_A-=A_{ij}\otimes_A-$.
	
	\item 2-morphisms between $F_{ij}$ and $F_{mn}$ are considered to be bimodule homomorphisms. For 2-morphisms to or from the identity, we consider the identity as tensoring with $A$, and take bimodule homomorphisms in this case (these are technically bimodule homomorphisms in $(A$-$A)$-$\biMod$, as $A$ is not necessarily a finite dimensional $A$-module).\end{itemize}

We will not be working directly with $\cC_A$ because the endomorphism hom-space of the identity is not in general necessarily finite dimensional, and indeed is not necessarily of countable dimension. For example, let $\Gamma$ be a quiver on $\bbZ$ where there is one arrow $i\to i$ for each $i\in\bbZ$ and no other arrows. Define $A=\bbk\Gamma/\bbk\Gamma_2$. Then the image of each $e_i$ under a bimodule homomorphism $\varphi:A\to A$ is independent of the image of any other $e_i$, and as $\dim e_iAe_i=2$, this implies that $\dim\End_{\cC_A}(\bbon_\tti)\geq 2^{|\bbZ|}$. While this is not a connected algebra, we will also show later that the same inequality holds for $\End_{\cC_{A_{\zig}}}(\bbon_\tti)$, though the reasoning is more complicated.\par

To fix this, we introduce a generalisation of 2-category $\cC_{A,X}$ first defined in \cite{mazorchuk2016endomorphisms} Section 4.5:

\begin{defn}\label{WFCAX} Let $\cC_A$ be as defined above. Let $Z$ denote the subalgebra of $\End_{\cC_A}(\bbon_*)$ generated by $\id_{\bbon_*}$ and by any 2-morphism that factors over a non-identity 1-morphism. Let $X$ denote a local subalgebra of $\End_{\cC_A}(\bbon_*)$ of finite or countably infinite dimension containing $Z$. The 2-category $\cC_{A,X}$ is defined to have the same objects, 1-morphisms and 2-morphisms as $\cC_A$ with the exception that $\End_{\cC_{A,X}}(\bbon_*)=X$.\end{defn}

This definition makes sense, since $\Hom_{\cC_A}(\bbon_*,F_{ij})$ has finite dimension for any $i,j\in V_A$, and thus $Z$ is of at most countably infinite dimension.

\begin{psn}\label{CAXK-S} $\cC_{A,X}$ is a locally Krull-Schmidt 2-category.
	\begin{proof} $\cC_{A,X}$ is an locally additive, locally idempotent complete 2-category by definition. Since $F_{ij}$ is indecomposable and $\End_{\cC_{A,X}}(F_{ij})$ is finite dimensional, $\End_{\cC_{A,X}}(F_{ij})$ is local by a standard argument. We chose $\End_{\cC_{A,X}}(\bbon_*)$ to be local, completing the proof.\end{proof}\end{psn}

We will show that such $\cC_{A,X}$ exist for certain classes of path algebras.

\begin{psn}\label{bpAlgrad1} Let $A$ be a bound path algebra with underlying quiver $\Gamma_A$. Then $\rad A$ is generated as an ideal by the equivalence classes of all paths in $\Gamma_A$ of length at least 1.
	\begin{proof} Fix a basis $\ceB$ of $A$ such that $e_i\in\ceB$ for all $i\in V_A$, and let $E_{A}=\{e_i|i\in V_A\}$. Let $R_A$ be the ideal of $A$ generated by the equivalence classes of all paths in $\Gamma_A$ of length at least 1. We recall that $\rad(A)=\{r\in A|rS_i=0\forall i\in V_A\}$, where $S_i$ is the simple module $\bbk\hat{e_i}$ defined previously. Let $r\in \rad(A)$. If $r=me_i+\sum\limits_{p\in \ceB\setminus\{e_i\}} m_p p$ for $m,m_p\in\bbk$ and some $i\in V_A$, then $rS_i=mS_i=0$, and thus $m=0$. Therefore, since $A$ is a bound path algebra, it follows that $r\in R_A$ and thus $\rad(A)\subseteq R_A$.\par

Conversely, if $r\in R_A$, then since $A$ is a bound path algebra $r=\sum\limits_{p\in \ceB\setminus E_{A}} m_p p$ for $m_p\in \bbk$, and thus $rS_i=0$ for any $i\in V_A$. Therefore $r\in\rad(A)$ and $\rad(A)=R_A$ as required.\end{proof}\end{psn}

\begin{cor}\label{bpAlgradk} Let $A$ be a bound path algebra with underlying quiver $\Gamma_A$. Then $\rad^k A$ is generated as an ideal by the equivalence classes of all paths in $\Gamma_A$ of length $k$. In particular, there is some integer $m$ such that $\rad^m A=0$.
	\begin{proof} For the first claim, we proceed by induction. \autoref{bpAlgrad1} provides the base case. Assume that $\rad^{k-1} A$ is generated as an ideal by the equivalence classes of all paths of length $k-1$. It immediately follows that $\rad^{k} A$ contains the equivalence classes of every path of length $k$. Conversely, let $a\in\rad^k A$. Then $a=\sum\limits_{i=1}^m r_ib_i$ for some $m$, where $r_i\in\rad A$ and $b_i\in\rad^{k-1} A$ for all $i$. But then $a$ is a sum of elements of the form $xp_1yp_{k-1}z$, where $x,y,z\in A$, $p_1$ is (an equivalence class of) a path of length 1 and $p_{k-1}$ is (an equivalence class of) a path of length $k-1$. In particular, this summand can further be written as a sum of elements of the form $vp_k w$, where $p_k$ is (an equivalence class of) a path of length $k$. The first result follows.\par
		
		For the second statement, we note from the definition of $A$ that $A=\bbk\Gamma_A/I$ with $(\bbk\Gamma_A)_k\subseteq I$ for some finite $k$. That is, the equivalence class of every path of length at least $k$ is zero, and the statement follows.\end{proof}\end{cor}

\begin{psn}\label{BPAZLoc} Let $A$ be a bound path algebra and let $Z\subseteq\End_{(A\text{-}A)\text{-}\operatorname{biMod}}(A)$ be as defined in \autoref{WFCAX}. Then $Z$ is a local $\bbk$-algebra.
	\begin{proof} If $V_A$ is finite then this is the finitary case which has already been proved in \cite{mazorchuk2016endomorphisms} Section 4.5. Therefore assume $V_A=\bbZ$. We claim that the subspace $I$ generated by all bimodule endomorphisms that factor over some $A_{ij}$ is a maximal proper left ideal of $Z$. If $I$ is a proper left ideal, it is immediately maximal - if $J\supset I$, then there must be an element of the form $\id_A+v\in J$ for $v\in I$. But then $v\in I\Rightarrow v\in J$, and therefore $\id_A\in J$ and $J=Z$. That $I$ is an ideal is clear - by definition it is closed under addition, and composing an endomorphism that factors over some $A_{ij}$ with another endomorphism still results in an endomorphism that factors over the same $A_{ij}$, and hence $I$ is closed under composition with elements of $Z$. It remains to show that $I$ is proper. \par
		
		Assume for contradiction that $\id_A$ is a member of $I$. Then without loss of generality $\id_A=\sum\limits_{l=1}^s \tau_l\sigma_l$ for $\sigma_l: A\to A_{i_l j_l}$ and $\tau_l:A_{i_l j_l}\to A$ for $i_l,j_l\in\bbZ$ and $s$ finite. We can therefore construct the column morphism $\sigma=(\sigma_l)_{l=1,\dots,s}:A\to\bigoplus\limits_{l=1}^s A_{i_lj_l}$ and the row morphism $\tau=(\tau_l)_{l=1,\dots,s}:\bigoplus\limits_{l=1}^s A_{i_lj_l}\to A$ such that $\tau\sigma=\id_A$. In particular, this implies that $\sigma_l$ is a split monomorphism, and thus in particular a monomorphism. This is a contradiction since $A$ is infinite dimensional and $\bigoplus\limits_{l=1}^s A_{i_lj_l}$ is finite dimensional. Therefore $I$ is indeed proper.\par
		
		We now claim that it is the unique maximal left ideal. By a straightforward but tedious consideration of cases, it can be shown that for any $v\in I$, $\operatorname{Im}(v)\subseteq \rad(A)$. Applying \autoref{bpAlgradk}, it follows that $v$ is nilpotent. Let $J$ be a maximal proper left ideal of $Z$. If $J\neq I$, then it must contain an element of the form $\id_A+v$ for $v\in I$. But since $v$ is nilpotent (say with nilpotency degree $k$), $$(\sum\limits_{j=0}^{k-1} (-v)^j)(\id_A+v)=\id_A+(-1)^k v^k=\id_A\in J,$$ a contradiction. Hence $I$ is the unique maximal proper left ideal of $Z$ and we are done.\end{proof}\end{psn}

\begin{psn}\label{QAlgcCWWF} Assume that $A$-$\modd$ is a Frobenius category. Then $\cC_{A,X}$ is a (locally) wide weakly fiat 2-category.
	
	\begin{proof} If $V_A$ is finite, then the statement has been proved in \cite{mazorchuk2011cell} Lemma 45 and \cite{mazorchuk2015transitive} Section 4.1. Therefore assume $V_A=\bbZ$. Let $*$ denote the unique object in $\cC_{A,X}$. It is immediate from the definitions that $\cC_{A,X}(*,*)$ is additive and $\bbk$-linear. Since the isomorphism classes of $F_{ij}$ are in bijection with $\bbZ\times\bbZ$ and there is a single isomorphism class of identity 1-morphisms, there are countably many indecomposable 1-morphisms up to isomorphism. We proved in \autoref{CAXK-S} that $\cC_{A,X}$ is Krull-Schmidt, and thus $\cC_{A,X}$ is a locally wide finitary 2-category.\par
		
		To show that $\cC_{A,X}$ is weakly fiat it is sufficient to show that each $F_{ij}$ is part of an internal adjunction. Since $A$-$\modd$ is a Frobenius category, by \autoref{BPAProj} it follows that the dual of any $Ae_i$ in $A$-$\modd$ is isomorphic to $Ae_{\sigma(i)}$ for some permutation $\sigma$ of $\bbZ$, which we refer to as the Nakayama bijection. We claim that the right internal adjoint of $F_{ij}$ is $F_{\sigma(j)i}$. To show this, we mirror the proof of \cite{mazorchuk2011cell} Lemma 45, adjusted to our setup. Given some $M\in A$-$\modd$, we have that
		\begin{align*} \Hom_{A\text{-}\modd}(Ae_i\otimes_\bbk e_jA, M) &\cong\Hom_{\bbk\text{-}\modd}(e_jA, \Hom_{A\text{-}\modd}(Ae_i,M))\\
		&\cong\Hom_{\bbk\text{-}\modd}(e_jA,e_iM)\\
		&\cong\Hom_{\bbk\text{-}\modd}(e_jA, e_iA\otimes_A M)\\
		&\cong\Hom_{\bbk\text{-}\modd}(e_jA,\bbk)\otimes_\bbk e_iA\otimes_A M\\
		&\cong(e_jA)^*\otimes_\bbk e_iA\otimes_A M.\end{align*}
		
		But as noted, $(e_jA)^*\cong Ae_{\sigma(j)}$, giving the claim.
\end{proof}\end{psn}

By \cite{dubsky2017koszulity} Remark 4, $A_{\zig}$-$\modd$ is a Frobenius category. We consider the 2-category $\cC_{A_{zig}}$ explicitly by examining the bimodule homomorphisms. It is straightforward to show the following: \begin{itemize}
	\item $\Hom_{\cC_{A_{\zig}}}(F_{ij},F_{mn})=0$ whenever $|i-m|\geq2$ or $|j-n|\geq 2$. For the case $|i-m|=|j-n|=1$ the hom-space is 1-dimensional, if $|i-m|+|j-n|=1$ it is 2-dimensional and if $|i-m|=|j-n|=0$ then the hom-space is 4-dimensional.
	
	\item If $|i-j|>2$, at least one of $e_kA_{\zig}e_i$ and $e_jA_{\zig}e_k$ is zero, and hence $\Hom_{\cC_{A_{\zig}}}(\bbon_*, F_{ij})=\Hom_{\cC_{A_{\zig}}}(\bbon_*,A_{\zig}e_i\otimes_\bbk e_jA_{\zig})=0$. If $|i-j|=2$ the hom-space is 1-dimensional, if $|i-j|=1$ the hom-space is 4-dimensional and if $i=j$ the hom-space is 6-dimensional.
	\item If $|i-j|\geq 2$, $\Hom_{\cC_{A_{\zig}}}(F_{ij},\bbon_\tti)=0$. If $|i-j|=1$ the hom-space is 1-dimensional, and if $i=j$ the hom-space is 2-dimensional.\end{itemize}

It remains to examine $\End_{(A_{\zig}\text{-}A_{\zig})\text{-}\biMod}(A_{\zig})$ and $Z_{\zig}$. If we have a bimodule endomorphism $\varphi$ such that $\varphi(e_i)=e_i$ for some $i$, then $\varphi(a_i)=a_i\varphi(e_i)=a_i$, but $\varphi(a_i)=\varphi(e_{i+1})a_i$ which implies that $\varphi(e_{i+1})=e_{i+1}$. A similar argument applies for $b_i$ and $e_{i-1}$, and therefore by bidirectional induction we derive that $\varphi(e_i)=e_i$ for all $i\in\bbZ$, and hence $\varphi=\id_A$.\par

If $\varphi(e_i)=b_ia_i$, then $\varphi(a_i)=0$ and thus $\varphi(e_{i+1})a_i=0$, and thus $\varphi(e_{i+1})$ is either (a scalar multiple of) $b_{i+1}a_{i+1}$ or $0$. However, this choice can be made freely, and it follows that a basis of $\End_{(A_{\zig}\text{-}A_{\zig})\text{-}\biMod}(A_{\zig})$ consists of the identity and of homomorphisms of the form $\varphi_I$, where $I\subseteq\bbZ$ and $\varphi_I(e_i)=b_ia_i$ if $i\in I$, and $0$ otherwise. The set of these $\varphi_I$ is thus in bijection with the powerset of $\bbZ$, and hence $\End_{\cC_{A_{\zig}}}(\bbon_\tti)$ has uncountable dimension.\par

Regarding $Z_{\zig}$, we can write $\varphi_{\{i\}}$ as $\sigma_i\tau_i$, where $\tau_i:A_{\zig}\to A_{\zig}e_{i+1}\otimes_\bbk e_{i+1}A_{\zig}$ is given by $\tau_i(e_j)=\delta_{ij} b_i\otimes a_i$, and $\sigma_i:A_{\zig}e_{i+1}\otimes_\bbk e_{i+1}A_{\zig}\to A_{\zig}$ given by $\sigma_i(e_{i+1}\otimes e_{i+1})=e_{i+1}$. Consequently, $Z_{\zig}$ is generated by $\id_A$ and by those $\varphi_I$ where $I$ is a finite subset of $\bbZ$. By \autoref{BPAZLoc}, $Z$ is a local algebra and therefore $\cC_{A_{\zig},Z}$ is a locally wide weakly fiat 2-category.

In general, $$A_{ij}\otimes_A A_{kl}\cong Ae_i\otimes_\bbk e_jAe_k\otimes_\bbk e_lA\cong (A_{il})^{\oplus\dim e_jAe_k}.$$ It follows that $F_{ij}\circ F_{kl}\cong F_{il}^{\oplus\dim e_jAe_k}$. In particular, the $\cL$-cells of $\cC_A$ are of the form $\ceL_j=\{F_{ij}|i\in\bbZ\}$ (and $\ceL_*=\{\bbon_*\}$) while the $\cR$-cells are of the form $\ceR_i=\{F_{ij}|j\in\bbZ\}$ (and $\ceR_*=\{\bbon_*\}$). It thus follows that $F_{ij}\sim_{\cL\circ\cR} F_{mn}$ for any $i,j,m,n\in\bbZ$, and hence there are two $\cD$-cells: $\cD_{*}=\{\bbon_*\}$ and $\cD_\bbZ=\{F_{ij}|i,j\in\bbZ\}$. But since it is clear that $\bbon_* >_{\cJ} F_{ij}$ for any $i,j\in\bbZ$, it follows that, on $\cC_{A,X}$, the $\cJ$ partial order agrees with the $\cD$ partial order, and the $\cJ$-cells are precisely the $\cD$-cells. Further, since $\ceL_j\cap\ceR_i=\{F_{ij}\}$, both $\cJ$-cells are strongly regular.\par

This allows us to construct the cell 2-representations corresponding to $\cC_{A,X}$. Choose some $j\in\bbZ$ and consider the $\cL$-cell $\ceL_j$. Then $$\add\{FX|F\in\cC, X\in\ceL_j\}=\add\{F_{mn}F_{ij}|m,n,i\in\bbZ\}=\add\ceL_j.$$ We denote this 2-representation of $\cC_A$ by $\mbN_j$. To recall some notation, we define the bimodule homomorphism $\varphi_{a,b}:A_{ij}\to A_{kl}$ by $\varphi_{a,b}(e_i\otimes e_j)=a\otimes b$. These $\varphi_{a,b}$ span $\Hom_{(A\text{-}A)\text{-}\operatorname{bimod}}(A_{ij},A_{kl})$. We have the following useful result:

\begin{lem}\label{WFCellRad} The maximal ideal $\ceI$ of $\mbN_j$ that is $\cC_{A,X}$-stable and does not contain any identity morphism for non-zero objects is generated as a collection of $\bbk$-vector spaces by the set $\{\varphi_{a,b}|b\in\operatorname{rad} e_jAe_j\}$.
	\begin{proof} The proof of this result generalises mutatis mutandis from the proof given for \cite{macpherson2020extension} Proposition 3.10.\end{proof}\end{lem}

In the case of $\cC_{A_{\zig}}$, $\rad e_jAe_j$ is a 1-dimensional vector space spanned by $a_jb_j$.\par

We denote the cell 2-representation corresponding to $\ceL_j$ by $\mbC_j$. In this case, we will show that we can give a stronger result that \autoref{ProCoalgExist}. Specifically, we will show that for an object $S\in\mbC_j(\tti)$, $[S,S]$ is a coalgebra 1-morphism in $\cC_{A,X}$ (or more precisely, its image under the forgetful functor $\Forg_S:\comod([S,S])\to\Pro(\widehat{\cC_{A,X}})$ lives in the image of $\cC_{A,X}$ under the canonical injection 2-functor). Let $\ceC_{A,X}=\cC_{A,X}(*,*)$.

\begin{psn}\label{eviotainC} The functor $\ev_j:\ceC_{A,X}\to\mbC_j(*)$ given by $\ev_j(F)=FF_{jj}$ and $\ev_j(\alpha)=\alpha\circ_H\id_{F_{jj}}$ is right adjoint to the functor $\Forg:\mbC_j(*)\to\ceC_{A,X}$ given by $\Forg(F)=F$ and $\Forg(\alpha)=\alpha$.
	
	\begin{proof} We will prove the adjunction by constructing the unit and counit adjunctions. By injection-projection arguments, it is sufficient to consider the components of the counit and the unit on indecomposable 1-morphisms/indecomposable objects. By a similar argument to the working in the proof of \cite{macpherson2020extension} Proposition 3.10, $\Hom_{\cC_{A,X}}(F_{ij},\bbon_*)$ consists of homomorphisms of the form $\varphi_a:A_{ij}\to A$ where $\varphi_a(e_i\otimes e_j)=a\in e_iAe_j$. On the other hand, $\Hom_{\mbC_j}(F_{ij},F_{jj})$ consists of linear combinations of homomorphisms of the form $\varphi_{a,e_j}:A_{ij}\to A_{jj}$ where $\varphi_{a,e_j}(e_i\otimes e_j)=a\otimes e_j$ for $a\in e_iAe_j$. Similarly, the morphism space $\Hom_{\mbC_j}(F_{ij},F_{kl}F_{jj})$ consists of linear combinations of morphisms of the form $\varphi_{a,b,e_j}$ where $\varphi_{a,b,e_j}(e_i\otimes e_j)=a\otimes b\otimes e_j$ for $a\in e_iAe_k$ and $b\in e_kAe_j$. In addition, again by a similar argument to the working in the proof of \cite{macpherson2020extension} Proposition 3.10, the morphism space $\Hom_{\cC_{A,X}}(F_{ij},F_{kl})$ consists of linear combinations of morphisms of the form $\varphi_{a,b}:A_{ij}\to A_{kl}$ where $\varphi_{a,b}(e_i\otimes e_j)=a\otimes b$ with $a\in e_iAe_k$ and $b\in e_lAe_j$.\par
		
		We define the unit morphisms $\eta_{F_{ij}}\in\Hom_{\mbC_j(*)}(F_{ij},F_{ij}F_{jj})$ as $\eta_{F_{kl}}=\varphi_{e_i,e_j,e_j}$. For the counit $\epsilon$, we first define $\epsilon_{\bbon_*}\in\Hom_{\cC_{A,X}}(F_{jj},\bbon_*)$ as $\epsilon_{\bbon_*}=\varphi_{e_j}$. We then define $\epsilon_{F_{kl}}\in\Hom_{\cC_{A,X}}(F_{kl}F_{jj},F_{kl})$ as $\epsilon_{F_{kl}}=\id_{F_{kl}}\circ_H\epsilon_{\bbon_*}$. It is straightforward to show that these satisfy the unit/counit axioms, and the results follows. \end{proof}\end{psn}

\begin{cor}\label{basecoalgwidef} Let $\mbC_j$ be a cell 2-representation of $\cC_{A,X}$. Then there exists some object $S\in\coprod\limits_{\tti\in\cC_{A,X}}(\mbC_j(\tti))$ such that the restriction of $\operatorname{Forg}_S$ to $[S,\coprod\limits_{\tti\in\cC}(\mbC_j(\tti))]$ factors over $\coprod\limits_{\ttj\in\cC_{A,X}}\cC_{A,X}(\tti,\ttj)$.
	\begin{proof} By \autoref{ProComodEquiv} for any $S\in\ceC_j$ there is an equivalence of 2-representations between $\mbC_j$ and $[S,\mbC_j]$, where $[S,-]$ is the right adjoint of $\Pro(\widehat{\ev_S})$ as defined previously. By \autoref{ProEvnN}, $\Pro(\widehat{\ev_S})$ is evaluation at $S$. If we set $S=F_{jj}$, then it follows by \autoref{eviotainC} that $[S,\ceC_j]=\Forg(\ceC_j)\subseteq\ceC_{A,X}$ as required.\end{proof}\end{cor}

We can say more. To begin, we present the generalisation of \cite{mackaay2016simple} Corollary 4.10 to the locally wide fiat case:

\begin{psn} Let $\cC$ be a locally wide weakly fiat 2-category, and let $\tti\in\cC$. Denoting by $\cC_\tti$ the endomorphism 2-category of $\tti$ in $\cC$, there is a bijection between equivalence classes of simple transitive 2-representations of $\cC_\tti$ and equivalence classes of simple transitive 2-representations of $\cC$ that have a non-trivial value at $\tti$.
	\begin{proof} For a simple transitive 2-representation $\mbM$ of $\cC$ and some object $S\in\ceM$, the (injective version of the) proof given in \cite{mackaay2016simple} assumes only the existence of the equivalence of 2-representations between $\mbM$ and $[S,\mbM]$, given here in \autoref{ProComodEquiv}. Beyond that we only need that for $S\in\mbM(\tti)$, $[S,S]$ lives in $\Pro(\widehat{\cC_{A,X}})(\tti,\tti)$, which remains true by definition.\end{proof}\end{psn}

We can now generalise \cite{macpherson2020extension} Lemma 3.11 to our setup:

\begin{thm}\label{stcellCAX} Any simple transitive 2-representation of $\cC_{A,X}$ is equivalent to a cell 2-representation.
	\begin{proof} We consider a larger 2-category $\cC_{A\times\bbk,X\times\bbk}$, defined similarly to $\cC_{*,\bbk}$ as in the proof of \cite{macpherson2020extension} Lemma 3.11. The endomorphism 2-category $\cC_\bbk$ of $*_\bbk$ is identical to that found in the referenced proof, and thus in particular all simple transitive 2-representations of it are equivalent to the cell 2-representation on it. Finally, any 1-morphism $Ae_i\otimes_\bbk e_jA$ of $\cC_{A,X}$ still factors over $*_\bbk$, and hence the rest of the proof of \cite{macpherson2020extension} Lemma 3.11 generalises without issue. \end{proof}\end{thm}

\begin{cor}\label{WFPathNice} Let $\mbM$ be a simple transitive 2-representation of $\cC_{A,X}$. Then there exists some object $S\in\ceM$ such that the restriction of $\operatorname{Forg}_S$ to $[S,\ceM]$ factors over $\ceC$.
	\begin{proof} This is a direct consequence of combining \autoref{basecoalgwidef} and \autoref{stcellCAX}\end{proof}\end{cor}

\section{Application: Soergel Bimodules}\label{WFSoerg}

We give a second application of the theory by demonstrating that the 2-category associated to a collection of Soergel bimodules is in fact a (locally) wide finitary 2-category.

\subsection{Soergel Bimodules: the Definitions}

We begin by defining Soergel bimodules, which were originally defined in \cite{soergel1992combinatorics}, though we draw our definitions here more from the summary paper \cite{williamson2011singular}, with some alterations for our specific case.

\begin{defn} Given $\bbZ$-graded algebras $A$ and $B$, we denote by $(A\text{-}B)\text{-}\tensor*[_\bbZ]{\biMod}{_{\bbZ,0}}$ the category whose objects are graded $(A$-$B)$-bimodules and whose morphisms are homogeneous graded bimodule homomorphisms of degree zero.\end{defn}

\begin{defn} A \emph{(finite) Coxeter matrix} $M$ is a symmetric square matrix with entries in $\bbZ^+\cup\{\infty\}$ where the diagonal entries are $1$ and the non-diagonal entries are at least $2$.\end{defn}

\begin{defn} A \emph{Coxeter system} is a pair $(W,S)$ where $W$ is a group and $S$ is a finite subset such there is a presentation of $W$ of the form $$\langle s\in S|(sr)^{m_{sr}}=e\text{ whenever }m_{ij}\text{ is finite.}\rangle$$ where $M=(m_{sr})_{s,r\in S}$ is a Coxeter matrix. The elements of $S$ are called \emph{simple reflections} and any element $w\in W$ that is conjugate to some $s\in S$ is called a \emph{reflection}.\end{defn}

\begin{defn} Let $(W,S)$ be a Coxeter system and let $x\in W$. An \emph{expression} of $x$ is a tuple $(s_1,s_2,\dots,s_n)\in S^n$ for some finite $n$ such that $x=s_1s_2\dots s_n$. This expression is \emph{reduced} if $n$ is minimal, and in this case we call $n$ the \emph{length} of $x$, denoted $l(x)$.\end{defn}

\begin{defn} Given a subset $I\subset S$, we denote by $W_I$ the subgroup of $W$ generated by $I$. This is generally called the \emph{parabolic subgroup} associated to $I$. By construction, $(W_I, I)$ is a Coxeter system. If $W_I$ is a finite subgroup, we say that $I$ is \emph{subgroup-finite}. If $I$ is subgroup-finite, we denote by $w_I$ the unique longest element of $W_I$.\end{defn}

Note that \cite{williamson2011singular} calls subgroup-finite subsets `finitary' subsets, but this presents obvious confusion issues for this paper, leading to the alternate nomenclature.

\begin{defn} For subsets $I,J\subseteq S$ let $W_I\setminus W/W_J$ denote the set of double cosets of $W$ with respect to $W_I$ and $W_J$. For a double coset $p\in W_I\setminus W/W_J$, we let $p_-$ denote the unique element of minimal length in $p$.\end{defn}

\begin{defn} Let $V$ be a ($\bbk$-)representation of $W$. For a subset $X\subset W$, we let $V^X$ denote the subspace of $V$ invariant under every element of $X$. We generally notate $V^{\{w\}}= V^w$ for simplicity.\end{defn}

\begin{defn} Let $V$ be a finite dimensional representation of $W$. We say $V$ is a \emph{reflection faithful} representation of $W$ if:\begin{itemize}
		\item The representation is faithful.
		\item $\codim V^w=1$ if and only if $w$ is a reflection.
	\end{itemize}
\end{defn}

By \cite{soergel2007kazhdan} Proposition 2.1, any Coxeter system has a reflection faithful representation. Let $R=S(V^*)$ be the symmetric algebra on $V^*$, graded such that $V^*$ is in degree 2. There is a natural action of $W$ on $R$ given by $w(f(v))=f(w^{-1}v)$ for any $w\in W$, $f\in R$ and $v\in V$. Thus given any subset $I\subset S$ or element $w\in W$ we can define $R^w$ and $R^I$ to be the subalgebras of $R$ invariant under $w$ and $W_I$ respectively. We note that, since $\bbk$ is of characteristic 0, $R$ is a graded-free $R^I$-module for any $I\subseteq S$. We recall the following result from \cite{williamson2011singular}:

\begin{psn}[\cite{williamson2011singular} Lemma 4.1.3] For $I\subseteq J$ subgroup-finite subsets of $S$, $R^I$ is a finitely generated graded free $R^J$-module. In addition, $$\Hom_{R^J\text{-}\tensor*[_\bbZ]{\Mod}{}}(R^I\llbracket l(w_I)-l(w_J)\rrbracket, R^J)\cong R^I\llbracket l(w_I)-l(w_J)\rrbracket.$$\end{psn}

\begin{defn} Let $I,J\subseteq S$ and let $p\in W_I\setminus W/W_J$. Set $K=I\cap p_-J{p_-}^{-1}$. The \emph{standard module indexed by $(I,p,J)$}, which we notate by $\tensor*[^I]{R}{_p^J}$, is an object of $(R^I$-$R^J)$-$\tensor*[_\bbZ]{\biMod}{_{\bbZ,0}}$ that as a ring is equal to $R^W$. The bimodule actions are given by:\begin{itemize}
		\item $r\cdot m=rm$ for $r\in R^I$ and $m\in \tensor*[^I]{R}{_p^J}$;
		\item $m\cdot r=m(p_-r)$ for $r\in r^J$ and $m\in\tensor*[^I]{R}{_p^J}$.
\end{itemize}\end{defn}
Of primary interest for us is the case where the double coset $p$ contains the identity $e$. In this case, we drop the $p$ from the notation and simply write $\tensor*[^I]{R}{^J}$.

\begin{defn} Let $I,J, K\subseteq S$ be subgroup-finite sets. We define three functors based on $J$ and $K$, all (following \cite{williamson2011singular}) notated as $\tensor*[^J]{\vartheta}{^K}$ based on the following conditions:\begin{itemize}
		\item If $J\subset K$ then we have $$\tensor*[^J]{\vartheta}{^K}:(R^J\text{-}R^I)\text{-} \tensor*[_\bbZ]{\biMod}{_{\bbZ,0}}\to(R^K\text{-}R^I)\text{-} \tensor*[_\bbZ]{\biMod}{_{\bbZ,0}}$$ given on objects by $\tensor*[^J]{\vartheta}{^K}(M)=\tensor*[_{R^K}]{M}{}\llbracket l(w_J)-l(w_K)\rrbracket$ with the obvious action on morphisms.
		\item If $K\subset J$ then we have $$\tensor*[^J]{\vartheta}{^K}:(R^J\text{-}R^I)\text{-}\tensor*[_\bbZ]{\biMod}{_{\bbZ,0}}\to(R^K\text{-}R^I)\text{-}\tensor*[_\bbZ]{\biMod}{_{\bbZ,0}}$$ given on objects by $\tensor*[^J]{\vartheta}{^K}(M)=R^K\otimes_{R^J}M$ and given on morphisms by $\tensor*[^J]{\vartheta}{^K}(f)=\id_{R^K}\otimes f$.
		\item If $K=J$, we set $\tensor*[^J]{\vartheta}{^K}$ to be the identity functor on $(R^J\text{-}R^I)\text{-} \tensor*[_\bbZ]{\biMod}{_{\bbZ,0}}$.
	\end{itemize}
	
	Since the definition of $\tensor*[^J]{\vartheta}{^K}$ has mutually exclusive components, there is never any ambiguity in its use, and this multi-part definition allows for more compact definitions later.\end{defn}

This allows us to define the bicategory of singular Soergel bimodules:

\begin{defn} The \emph{bicategory of singular Soergel bimodules} $\cB_b=\cB_{(W,S),b}$ has objects enumerated by the subgroup-finite subsets $I$ of $S$. We abuse notation by also referring to the object associated to $I\subseteq S$ as $I$. The hom-categories $\cB(I,J)$ are the smallest additive subcategories of $(R^J\text{-}R^I)\text{-}\tensor*[_\bbZ]{\biMod}{_{\bbZ,0}}$ subject to the following:\begin{enumerate}
		\item $\cB_b(I,J)$ is closed in $(R^J\text{-}R^I)\text{-}\tensor*[_\bbZ]{\biMod}{_{\bbZ,0}}$ under taking isomorphisms, direct summands and $\bbZ$-grading shifts.
		\item $\tensor*[^I]{R}{^I}$ is a 1-morphism in $\cB_b(I,I)$ for all objects $I\in\cB_b$.
		\item If $B\in\cB_b(I,J)$, then $\tensor*[^J]{\vartheta}{^K}(B)\in\cB_b(I,K)$ whenever this is defined (i.e. when $K$ is subgroup-finite with either $K\subseteq J$ or $J\subseteq K$).
	\end{enumerate}
	We take composition of 1-morphisms and horizontal composition of 2-morphisms to be tensor products over then common algebra (i.e. for $B\in\cB(I,J)$ and $C\in\cB(J,K)$, $C\circ B=C\otimes_{R^J} B\in\cB(I,K)$). We note that $\bbon_I=R^I$.\end{defn}

As noted in Section 1 of \cite{williamson2011singular}, this is equivalent to setting $\cB_b(I,J)$ as the smallest full additive sub-category of $(R^J\text{-}R^I)$-$\tensor*[_\bbZ]{\biMod}{_{\bbZ,0}}$ containing all objects isomorphic to direct summands of shifts of objects of the form $$R^{I_1}\otimes_{R^{J_1}} R^{I_2}\otimes_{R^{J_2}}\dots\otimes_{R^{J_{n-1}}} R^{I_n}$$ with $I=I_1\subset J_1\supset I_2\subset\dots\subset J_{n-1}\supset I_n=J$ subgroup-finite subsets.

\begin{defn} The \emph{2-category of singular Soergel bimodules} $\cB=\cB_{(WS)}$ is defined as a 2-category biequivalent to $\cB_b$.\end{defn}

The endomorphism sub-2-category $\cB_\emptyset$ of the object $\emptyset$ is referred to as the 2-category of \emph{Soergel bimodules}.

\subsection{Soergel Bimodules: the Structure}

We now demonstrate that the 2-category of singular Soergel bimodules is in fact a locally wide fiat 2-category. We first give a special case of \cite{elias2016soergel} Lemma 6.24, adapted to the language used in this thesis:

\begin{lem} The hom-category $\cB(\emptyset, \emptyset)$ is a Krull-Schmidt category.\end{lem}

The proof of \cite{elias2016soergel} Lemma 6.24 adapts without issue to any other hom-category in $\cB$, giving the following lemma:

\begin{lem}\label{SoergKS} For any objects $I$ and $J$ of $\cB$, $\cB(I,J)$ is a Krull-Schmidt category.\end{lem}

\begin{lem}\label{SSoergWFin} The 2-category $\cB$ is a locally wide finitary 2-category.
	\begin{proof} Since we took the Coxeter system $(W,S)$ such that $S$ is a finite set, $\cB$ has finitely many objects and thus certainly at most countably many objects. Since for any subgroup-finite $I$ and $J$ $\cB_b(I,J)$ is a sub-2-category of $(R^I\text{-}R^J)$-$\tensor*[_\bbZ]{\biMod}{_{\bbZ,0}}$ it is certainly additive and $\bbk$-linear with countable dimension hom-spaces of 2-morphisms, and hence so is $\cB(I,J)$.\par
		
		There are only countably many $(R^I\text{-}R^J)$-bimodules of the form  $$R^{I_1}\otimes_{R^{J_1}} R^{I_2}\otimes_{R^{J_2}}\dots\otimes_{R^{J_{n-1}}} R^{I_n}$$ for $I=I_1\subset J_1\supset I_2\subset\dots\subset J_{n-1}\supset I_n=J$ subgroup-finite subsets, and hence there are only countably many grade-shifts of these. Since each of these has only finitely many indecomposable direct summands, it follows that $\cB_b(I,J)$ (and hence $\cB(I,J)$) has countably many isomorphism classes of indecomposable 1-morphisms. By construction the identity 1-morphisms are indecomposable, and finally by \autoref{SoergKS} the hom-categories are Krull-Schmidt. This completes the proof.\end{proof}\end{lem}

In fact, we can say more:

\begin{lem}\label{SSoergWFiat} The 2-category $\cB$ is a locally wide fiat 2-category.
	\begin{proof} It is a consequence of \cite{soergel2007kazhdan} Proposition 5.10 that $R\otimes_{R^s} R\otimes_R-\in\cB(\emptyset, \emptyset)$ is self-adjoint for any $s\in S$, and it thus follows that the endomorphism 2-category of $\emptyset$ is locally wide fiat. That the auto-involution extends to the whole of $\cB$ comes from using a straightforward generalisation of \cite{soergel2007kazhdan} Proposition 5.10 to the singular Soergel bimodule setup.\end{proof}\end{lem}

This is a much wider array of Soergel bimodule 2-categories that can be studied than under simply finitary 2-representation theory. To give a clearer view of this, we will relate the Soergel bimodules to the Coxeter-Dynkin diagrams associated to the Coxeter groups. We will not discuss Coxeter-Dynkin diagrams in detail (see \cite{grove1985classification} Chapter 5 for a detailed study), but briefly they are graphs whose edges are labelled with positive integers which fully classify Coxeter groups.\par

Original finitary 2-representation theory can only cover finite Coxeter groups. These are classified by three infinite families of Coxeter-Dynkin diagrams (called $A_n$, $B_n$ and $D_n$), as well as a small finite set of exceptional diagrams. The wide finitary theory detailed herein instead applies to any Coxeter-Dynkin diagram with finitely many nodes, including not only the affine Coxeter-Dynkin diagrams, but also a wide array of wild Coxeter-Dynkin diagrams.

\begin{thm} Let $\cB=\cB_{(W,S)}$ be the 2-category of singular Soergel bimodules associated to a Coxeter system $(W,S)$ with $S$ a finite set. Let $\mbM$ be a transitive 2-representation of $\cB$. Then there exists a coalgebra 1-morphism $C$ in $\Pro(\widehat{\cB})$ such that $\mbM$ is equivalent as a 2-representation to a subcategory of the category of comodule 1-morphisms over $C$.
	\begin{proof} \autoref{SSoergWFiat} gives that $\cB$ is a locally wide fiat 2-category, and we can thus apply \autoref{ProComodEquiv} to get the result immediately.\end{proof}\end{thm}

\bibliography{C:/Users/jamac/Documents/Maths/Thesis/MasterBib}
\end{document}